\documentclass[12pt,reqno,twoside]{article}
\usepackage{amssymb,amsfonts,amsthm,amsmath}
\usepackage{graphicx}
\usepackage{color}
\usepackage{plain}
\usepackage{cite}
\usepackage[left=1 in,top=1 in,right=1 in,bottom=1 in]{geometry}

\numberwithin{equation}{section}

\def\beq{\begin{equation}}
\def\eeq{\end{equation}}
\def\beqs{\begin{equation*}}
\def\eeqs{\end{equation*}}
\def\dmax{\delta_{\rm max}}

\def\ca{C_1} \def\cb{C_2} \def\cc{C_3} \def\cd{C_4}
\def\ce{C_5} \def\cf{C_6} \def\cg{C_7} \def\ch{C_8}
\def\ci{C_9} 
\def\cl{C_{10}} 
\def\cj{C_{11}}
\def\ck{C_{12}}
\def\cm{C_{13}} \def\cn{C_{14}}
\def\co{C_{15}} \def\cp{C_{16}} \def\cq{C_{17}}
\def\cs{C_{18}}

\newtheorem{theorem}{Theorem}[section]
\newtheorem{lemma}[theorem]{Lemma}
\newtheorem{proposition}[theorem]{Proposition}
\newtheorem{corollary}[theorem]{Corollary}
\theoremstyle{definition}

\newtheorem{definition}[theorem]{Definition}

\def\A{{\cal A}}

\def\Z{{\bf Z}}
\def\C{{\bf C}}
\def\J{{\cal J}}

\def\words#1{\quad\hbox{#1}\quad}
\def\wwords#1{\qquad\hbox{#1}\qquad}

\numberwithin{equation}{section}

\title{
    A Relaxed Direct-insertion Downscaling Method\\
	For Discrete-in-time Data Assimilation
}

\author{Emine Celik$^1$ and Eric Olson$^2$}
\date{\today}

\begin{document}
\maketitle
\begin{center}
\textit{$^1$Department of Mathematics,
Sakarya University\\
54050 Sakarya, T\"{u}rkiye}\\
\textit{$^2$Department of Mathematics and Statistics,
University of Nevada, Reno\\
Reno, NV  89557, USA}\\
\medskip
Email addresses:  \texttt{eminecelik@sakarya.edu.tr, ejolson@unr.edu}
\end{center}

\begin{abstract}
\noindent
This paper improves the spectrally-filtered direct-insertion 
downscaling method for discrete-in-time data assimilation
by introducing a relaxation parameter that overcomes a constraint
on the observation frequency.
Numerical simulations demonstrate that taking the relaxation 
parameter proportional 
to the time between observations allows one to vary 
the observation frequency over a wide range while maintaining
convergence of the approximating solution to the reference solution.
Under the same assumptions we analytically prove
that taking the 
observation frequency to infinity 
results in the continuous-in-time nudging method.
\end{abstract}

\section{Introduction}\label{intro}

Our focus is that of a well-posed dissipative dynamical system
\begin{equation}\label{freerunning}
	{dU\over dt}={\cal F}(U)
\wwords{with} U(t_0)=U_0
\end{equation}
where the initial condition $U_0$ is unknown but for which a time
sequence of partial observations of $U(t)$ are available at times $t=t_n$.

Following \cite{AOT2014}, \cite{Titi2003}, related research, the 
references therein and in particular \cite{COT2019}, we consider the 
concrete setting of the incompressible 
two-dimensional Navier--Stokes equations.  
These equations provide an example of a well-posed
dissipative dynamical system of the 
form~\eqref{freerunning} that is
simpler than atmospheric or ocean models but with
similar nonlinear dynamics.
Other systems for which a growing body of analytic and 
numerical results are 
available include Rayleigh-B\'enard convection studied 
by Farhat, Lunasin and Titi in \cite{Farhat2016} (see also
Farhat, Glatt-Holtz, Martinez, McQuarrie
and Whitehead \cite{Farhat2019}) as well as the surface
quasi-geostrophic equations studied by Jolly, Martinez and Titi in
\cite{Jolly2017} (see also \cite{Jolly2019}) among others.
Although our computations and analysis focus on the 
two-dimensional Navier--Stokes equations,
we hope the resulting intuition and conclusions will apply to
operational models with practical applications.

We begin with the spectrally-filtered discrete-in-time 
downscaling data assimilation algorithm introduced 
in~\cite{COT2019} in the context of the incompressible 
two-dimensional Navier--Stokes equations.
This is a direct-insertion method that recovers unobserved lengthscales 
by inserting new observational measurements as they are available 
into the current estimate of the state.
The spectral filter helps ensure no high-resolution artifacts are present 
in the interpolated observational data that might damage the approximation
obtained by the numerical model as it is integrated forward in time.
From a theoretical point of view, the filtering is needed to obtain 
suitable estimates in the higher Sobolev norms used to show convergence of the 
approximating solution to the reference solution over time.

\begin{figure}[h!]
    \centerline{\begin{minipage}[b]{0.75\textwidth}
    \caption{\label{epaths}%
    The error $|U(t)-u(t)|$
    for a reference solution $U(t)$
    where the approximating solution $u(t)$ was computed 
	for different observation intervals $\delta$.
	Large $\delta$ is shown on the left; small on the right.}
    \end{minipage}}
    \centerline{
        \includegraphics[height=0.4\textwidth]{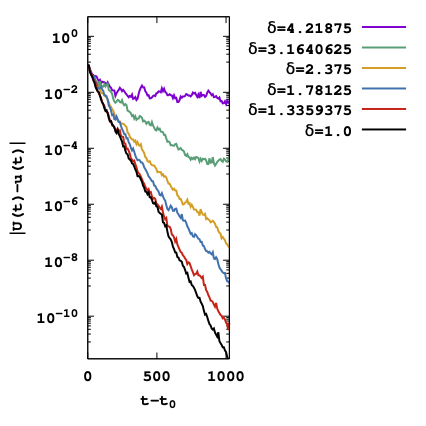}
        \hskip-14pt
        \includegraphics[height=0.4\textwidth]{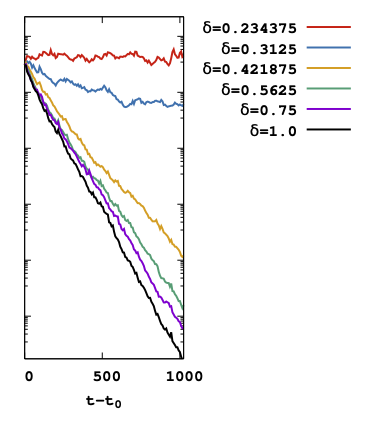}
    }
\end{figure}

Intuitively more frequent observations should make the predicted state 
obtained from 
data assimilation more accurate.
For example, 
in weather forecasting 
measuring the state of the atmosphere more frequently 
should---provided a reasonable algorithm is used---allow an 
approximating solution to track the reference solution more accurately.
However, although more frequent observations 
represent greater knowledge about the reference solution, 
Figure~\ref{epaths}
shows the spectrally-filtered discrete-in-time
data assimilation algorithm
introduced in \cite{COT2019} 
actually performs worse
when data is inserted more frequently into the model.

The possibility of this unreasonable behavior was 
already noticed 
in the theoretical analysis presented in \cite{COT2019} and
commented on as
\par\medskip
{\narrower\noindent\it
[Our analysis makes]
use of a minimum distance between $t_{n+1}$ and $t_n$ as well as the
maximum. Measurements need to be inserted frequently enough to overcome
the tendency for two nearby solutions to drift apart, while at the
same time the possible lack of orthogonality in our general interpolant
observables means measurements should not be inserted too frequently.
\par}
\medskip\noindent
The present research is motivated by numerical evidence that the above
constraint concerning the minimum observation frequency is physical
and not merely a limitation of our analytic techniques.
In particular, given the ideal situation of noise-free observations
of exact dynamics in which the approximating solution synchronizes with
the reference solution over time, our computations show that inserting 
measurements more 
frequently can worsen the quality of the approximation to the point where
it subsequently fails to synchronize.
To focus on the situation when subsequent observations become more frequent,
we assume $t_{n+1}-t_n\le\dmax$ for all $n$.

Recall the incompressible two-dimensional Navier-Stokes equations 
given by
\begin{equation}\label{nse0}
	\frac{\partial U}{\partial t}+(U\cdot \nabla)U
		-\nu \Delta U+\nabla p=f,\qquad
	\nabla \cdot U=0,
\end{equation}
where $U$ is the velocity of the fluid, $\nu$ is the kinematic viscosity,
$p$ is the pressure and $f$ is a time-independent body force applied to 
the fluid.
Consider the spatial domain $\Omega$ with $L$-periodic
boundary conditions for simplicity.

Recast~\eqref{nse0} in the functional setting
described by Constantin and Foias~\cite{Constantin1988} (see also
Robinson \cite{Robinson2001} or Temam \cite{Temam1977}) as follows.
Denote by ${\cal V}$ the set of all divergence-free $L$-periodic trigonometric
polynomials with zero spatial averages.  Let $V$ be the closure
of ${\cal V}$ in $H^1(\Omega,\mathbf{R}^2)$ and $H$ be the closure 
of ${\cal V}$
in $L^2(\Omega,\mathbf{R}^2)$.  Denote the dual of $V$ by $V^*$.

Let
$A\colon V\to V^*$ and $B\colon V\times V\to V^*$ be 
continuous extensions of
the operators given by $Au=-P_H\Delta u$ 
and $B(u,v)=P_H(u\cdot\nabla v)$ for $u,v\in{\cal V}$
where 
$P_H$ is the orthogonal projection 
of $L^2(\Omega)$ onto $H$.  
Let ${\cal D}(A)=\{ u\in V: Au\in H\}$ be the domain of $A$ into $H$.
Applying $P_H$ to both sides of \eqref{nse0}
then expresses the Navier--Stokes equations as
\beq
\label{nse1} \frac{dU}{dt}+\nu AU+B(U,U)=f
\wwords{where}
U(t_0)=U_0\in V.
\eeq
Note in the periodic case $A=-\Delta$ and that we have assumed $P_H f=f$.

The spaces $H$, $V$ and ${\cal D}(A)$ can be characterized 
in terms of Fourier series as $H=V_0$, $V=V_1$ and ${\cal D}(A)=V_2$ 
where 
$$
    V_\alpha=\bigg\{ \sum_{k\in\J} u_k e^{ik\cdot x}
        : 
        \sum_{k\in\J} |k|^{2\alpha}|u_k|^2<\infty,\quad
        k\cdot u_k=0
        \words{and}
        u_{-k}=\overline{u_k}
    \,\bigg\}.
$$
Here ${\cal J}=\big\{ \,2\pi n/L:n\in \Z^2\setminus \{0\}\big\}$ and
$u_k\in\C^2$ are the Fourier coefficients for the velocity field $u$.
The zero spatial average along with Parseval's identity implies the 
norms and inner products for the $V_\alpha$ spaces are equivalent to
$$
	\|u\|_\alpha = \Big(L^2\sum_{k\in\J} 
		|k|^{2\alpha}|u_k|^2\Big)^{1/2}
\wwords{and}
	(\!(u,v)\!)_\alpha= L^2\sum_{k\in\J} |k|^{2\alpha} u_k\overline{v_k}.
$$

Recall that $V_\alpha\subseteq V_\beta$ for $\beta\le\alpha$ with
corresponding Poincar\'e inequalities
\begin{equation}\label{poincare}
	\lambda_1^{\alpha-\beta}\|u\|_\beta^2\le \|u\|_\alpha^2
\wwords{for} u\in V_\alpha,
\end{equation}
where $\lambda_1=(2\pi/L)^2$.
Denote the norms corresponding 
to $H$, $V$ and ${\cal D}(A)$ by
$$
	|u|=\|u\|_0,\qquad
	\|u\|=\|u\|_1\wwords{and}
	|Au|=\|u\|_2
$$
and the inner products by
$$
	(u,v)=(\!(u,v)\!)_0,\qquad
	(\!(u,v)\!)=(\!(u,v)\!)_1\wwords{and}
	(Au,Av)=(\!(u,v)\!)_2.
$$

Along with the above functional notation the following {\it a priori\/} bounds
also appear in \cite{Constantin1988}, 
\cite{Robinson2001} and \cite{Temam1977}
for the incompressible Navier--Stokes equations:
\begin{theorem}\label{rhoapriori}
Let $U$ be a solution to \eqref{nse1} with $f\in V$ and initial 
condition $U_0\in {\cal D}(A)$.
There exist bounds
$\rho_\alpha$ and $\widetilde\rho_\alpha$
depending only on $\dmax$, $\nu$, $\|f\|$ and $U_0$ such that 
\begin{equation}\label{Urhov} 
	|A^{\alpha/2}U|\le \rho_\alpha
	\words{for}\alpha=0,1,2
\end{equation}
and
\begin{equation}\label{Uintbound}
	\Big(\int_{t_n}^{t_{n+1}}
		|A^{\alpha/2}U|^2
		dt\Big)^{1/2}	
	\le \widetilde \rho_\alpha
	\words{for}\alpha=1,2,3.
\end{equation}
Furthermore, \eqref{nse1} possesses a unique global attractor $\A$
bounded in ${\cal D}(A)$.  Moreover, if $U$ lies on that global attractor,
then the constants $\rho_\alpha$ and $\widetilde\rho_\alpha$ may 
be taken independent of $U_0$.
\end{theorem}
We pause to remark that 
a well-posed dissipative dynamical system of the form 
\eqref{freerunning} may now be obtained by 
taking ${\cal F}(U)=f-B(U,U)-\nu AU$.

Let $S$ be the solution operator such that 
$S(t)U_0=U(t)$, where $U(t)$ is the 
unique solution to~\eqref{nse1} with initial condition 
$U_0$ at time $t_0=0$.
Well-posedness along with the fact that 
the dynamics are autonomous implies
$S(t+\delta)=S(\delta)S(t)$ for $t\ge 0$ and $\delta\ge 0$.

The spectral filter appearing in \cite{COT2019} is given
for $u\in H$ as
\begin{equation}\label{filter}
	P_\lambda u = \sum_{0<|k|^2\le\lambda} u_k e^{ik\cdot x}
\words{where}
	u=\sum_{k\in\J} u_k e^{ik\cdot x}.
\end{equation}
Note the interpolation property
\begin{equation}\label{pinterp}
	|u-P_\lambda u|^2
	=L^2\sum_{|k|^2>\lambda} |u_k|^2
	\le L^2\sum_{|k|^2>\lambda} {|k|^2\over\lambda}|u_k|^2
	\le \lambda^{-1} \|u\|^2
\words{for} u\in V,
\end{equation}
and the smoothing properties
\begin{equation}\label{psmooth}
	\|P_\lambda u\|_\alpha^2=
		L^2\!\!\!\sum_{0<|k|^2\le\lambda} |k|^{2\alpha}|u_k|^2
		\le 
		\lambda^\alpha L^2\!\!\!\sum_{0<|k|^2\le\lambda} |u_k|^2
	\le \lambda^\alpha |u|^2
\words{for} u\in H.
\end{equation}
In particular note that
$\|P_\lambda u\|^2\le \lambda |u|^2$ and
$|AP_\lambda u|^2\le \lambda^2 |u|^2$.
Similarly $|AP_\lambda u|^2\le \lambda \|u\|^2$.

Assume the observations are interpolated back into the phase 
space of $U$ by a linear operator $I_h\colon V\to L^2(\Omega)$ that 
satisfies
\begin{equation}\label{interp}
	|U-P_H I_h U|^2\le c_1h^2 \|U\|^2
\words{for every} U\in V\words{and} h>0.
\end{equation}
Here $h$ is a length scale that reflects the resolution of the observations.
Note \eqref{interp} provides a bound on the relative error in the 
spirit of the Bramble--Hilbert Lemma and that there is no 
requirement on $I_h$ or $P_H I_h$ to be a projection.  
Linearity, however, is important.

The operators $I_h$ were introduced in \cite{AOT2014} 
as type-I interpolant observables 
and generally result from spatially-averaged observations,
for example, the determining volume elements of
Jones and Titi \cite{Jones1992}.
Note that observations which result from point measurements
are classified as type-II and rely on higher norms such 
as $|Au|$ to obtain inequalities similar to~\eqref{interp}.
Although the treatment of point 
measurements is out of scope for the present research,
we remark that
a physically motivated type-I interpolant involving locally 
averaged approximations of point measurements in space
is described in Appendix B of \cite{CO2022}.

Since $P_H I_h$ may not be an orthogonal projection,
lack of orthogonality complicates 
the use of $P_H I_h U$ for data assimilation.
Another
difficulty is the possible roughness of $I_h$.
Such roughness may be characterized as the case when
$P_H I_h U\not\in V$.
This difficulty was mediated in \cite{COT2019} by the 
use of $P_\lambda$ as a smoothing filter.
To this end set $J=P_\lambda P_H I_h$ where $P_\lambda$ is
the spectral filter defined by \eqref{filter}
and consider the smoothed observations given by $JU$.

This results in an interpolant that satisfies
\begin{equation}\label{Jinterp}
	|U-J U|^2=
	|U-P_\lambda U|^2+|P_\lambda (U-P_H I_h U)|^2
	\le
	(\lambda^{-1}+c_1h^2) \|U\|^2
\end{equation}
and
\begin{equation}\label{vcont}
	\|U-J U\|^2=
	\|U-P_\lambda U\|^2+\|P_\lambda (U-P_H I_h U)\|^2
	\le
	\big(1+\lambda c_1h^2\,\big) \|U\|^2,
\end{equation}
which has the continuity properties
\begin{equation}\label{Jweak}
	|J U|=|U-JU|+|U|\le (\lambda^{-1}+c_1 h^2)^{1/2}\|U\|+\lambda_1^{-1}\|U\|
		\le c_2 \|U\|
\end{equation}
and
\begin{equation}\label{Jsmooth}
    \|J U\|=\|U-JU\|+\|U\|
    \le (1+\lambda c_1h^2)^{1/2} \|U\| + \|U\|
	\le c_3 \|U\|
\end{equation}
where $c_2=(\lambda^{-1}+c_1 h^2)^{1/2}+\lambda_1^{-1/2}$ and
$c_3=(1+\lambda c_1h^2)^{1/2}+1$.

We now describe main theorem proved in \cite{COT2019} 
for the spectrally-filtered discrete-in-time data 
assimilation algorithm
whose improvement is the purpose of the 
present paper.

\begin{definition}\label{direct}
Given an increasing sequence $t_n$ of observation times,
the approximating solution $u(t)$ obtained by
{\it spectrally-filtered discrete-in-time data assimilation\/} is
$$
	u(t)=S(t-t_n)u_n\wwords{for} t\in[t_n,t_{n+1}),
$$
where
$$
\left\{\begin{aligned}
	u_0&=J U(t_0),\nonumber\\
	u_{n+1}&=
		(I-J) S(t_{n+1}-t_n)u_n+J U(t_{n+1})\label{nokappa}.
\end{aligned}\right.
$$
\end{definition}

We remark taking $J=P_\lambda$ leads to the
discrete-in-time data-assimilation method studied in 
\cite{Hayden2011} based on observation of the Fourier modes.
In this case the quality of the approximation $u(t)$ does not 
worsen when measurements are inserted more frequently in time.
Thus, lack of orthogonality and smoothness 
in $J$ appear to be important considerations.

When $t_n=t_0+\delta n$ for some period $\delta$,
the main result from \cite{COT2019} 
may be stated as

\begin{theorem}
\label{COTmain}
Let $U$ be a solution to the incompressible two-dimensional Navier--Stokes
equations and $u(t)$ be the process given by Definition~\ref{direct}.
Then, for every $\delta>0$ there exists $h>0$ and $\lambda>0$ depending
only on $c_1$, $f$ and $\nu$ such that
$$
	|u(t)-U(t)|\to 0\words{exponentially in time as}t\to\infty.
$$
Here $c_1$ is the constant appearing in \eqref{interp}.
\end{theorem}

In the presence of model error and noise
the synchronization 
given by Theorem~\ref{COTmain} 
as $t\to\infty$ 
would be only approximate.
While it seems reasonable that observing the reference solution less
frequently in time with larger $\delta$ can be compensated for by 
requiring higher resolution observations with smaller $h$, the
difficulty we address in this paper is that once the observation 
period $\delta$ is chosen,
the theory does not allow more frequent observations without 
again adjusting the resolution $h$ and filter parameter $\lambda$.

For simplicity assume as above that the time 
between subsequent observations is fixed to be $\delta$.
The case where $t_{n+1}-t_n=\delta_n$ with $\max_n \delta_n$ small
is interesting and commented on in the conclusions.

Before introducing a relaxation parameter into Definition~\ref{direct}
we first present additional numerical evidence 
to illustrate the constraints on the 
observation frequency suggested by the proof of 
Theorem~\ref{COTmain}.
Consider a fixed trajectory $U(t)$ on the global attractor of the 
two-dimensional Navier--Stokes equations obtained by a long-time 
integration of \eqref{nse1} forward in time.
Next, apply Definition \ref{direct} to
compute $u(t)$ for different values of $\delta$ 
while holding all other parameters and the trajectory
$U(t)$ fixed.

The evolutions of $|U(t)-u(t)|$ for the resulting simulations 
are depicted in Figure~\ref{epaths}.
Consider the value of $|U(T)-u(T)|$ at $T=t_0+1024$ for
the different values of $\delta$ and
the behavior of each curve leading up to that time.
The graph on the left suggests $|U(t)-u(t)|$ does not 
tend to zero as $t\to\infty$
for values of $\delta$ where $\delta>3$.
More surprisingly the graph on the right suggests $|U(t)-u(t)|$ 
does not tend to zero when $\delta<0.3$.
Thus,
taking $\delta$ either too small or too large results 
in failure of 
the approximating solution $u$ to synchronize with $U$.  
This is consistent with the proof of Theorem \ref{COTmain}.

\begin{figure}[h!]
    \centerline{\begin{minipage}[b]{0.75\textwidth}
    \caption{\label{epathsb}%
	The same calculation as Figure~\ref{epaths} 
	redone for a different reference 
	solution $U(t)$ lying 
    on the global attractor.}
    \end{minipage}}
    \centerline{
        \includegraphics[height=0.4\textwidth]{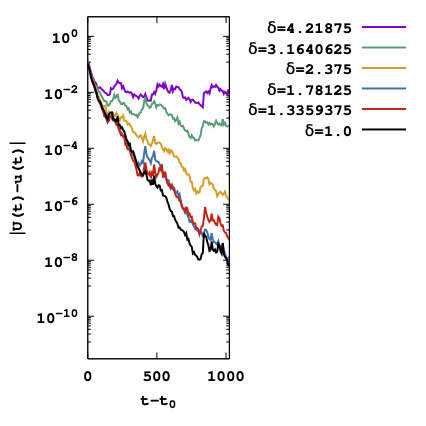}
        \hskip-14pt
        \includegraphics[height=0.4\textwidth]{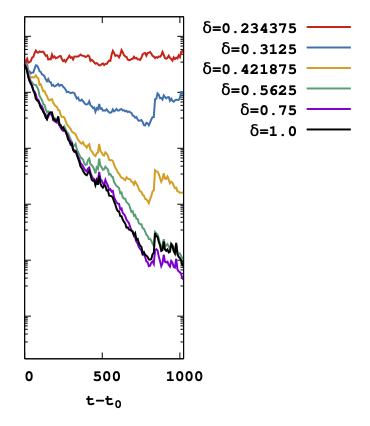}
    }
\end{figure}

Note the skill by which 
Definition \ref{direct} recovers the reference solution $U(t)$ 
depends on the specific trajectory being observed.
Some reference solutions may exhibit dynamical 
behaviors that are more difficult to recover from the observational
data while others are easier to recover.  
Figure~\ref{epathsb}
repeats the calculations of Figure~\ref{epaths} 
for a different reference solution on the global attractor.
While the size of $|U(T)-u(T)|$ at $T=t_0+1024$ 
is four decimal orders of magnitude larger 
in Figure~\ref{epathsb} compared to Figure~\ref{epaths},
the trends of the curves for large and small values of 
$\delta$ are similar.
In particular,
taking $\delta$ either too small or too large results 
in failure of 
the approximating solution $u$ to synchronize with $U$ 
in both cases while
$\delta\approx 1$ is a good
choice for the observation interval.

\begin{figure}[h!]
    \centerline{\begin{minipage}[b]{0.75\textwidth}
    \caption{\label{dbox}%
    The ensemble average and geometric mean compared to box plots
    (no outliers removed) of the error $|U(T)-u(T)|$ at time $T=t_0+1024$ 
    for $500$ simulations varying $\delta$.}
    \end{minipage}}
    \centerline{
        \includegraphics[height=0.4\textwidth]{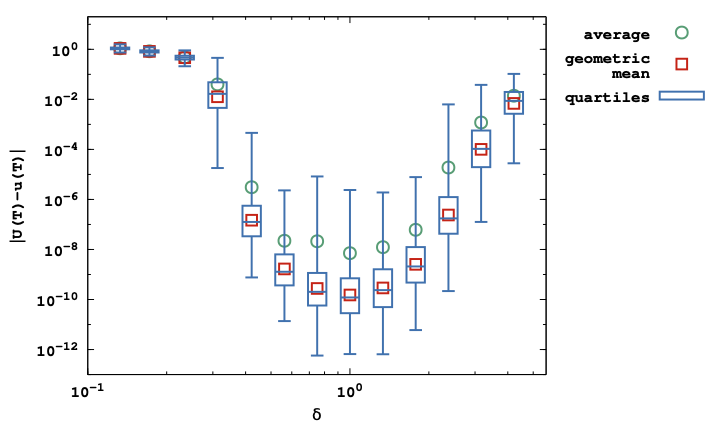}
    }   
\end{figure}

To further characterize the rate at which the approximating
solution converges to the reference solution for different
values of $\delta$,
the experiments of Figures~\ref{epaths} and~\ref{epathsb}
were repeated 
for an ensemble of
$500$ reference solutions randomly chosen on the global
attractor.
Let ${\cal E}$ be a fixed ensemble.
Compute the ensemble average error
for each $\delta$ as
$$\langle |U(T)-u(T)|\rangle
	={1\over |{\cal E}|}\sum_{U\in{\cal E}} |U(T)-u(T)|
$$
and the geometric mean as
$${\rm GM}\big[|U(T)-u(T)|\big]
	=\exp\Big({1\over |{\cal E}|}\sum_{U\in{\cal E}} \log|U(T)-u(T)|\Big).
$$
Here $|{\cal E}|$ is the number of solutions in the ensemble.

Figure \ref{dbox} reports the ensemble averages corresponding to
$\delta\in [0.1328125,4.21875]$ along with
the median and some additional statistics.
The whiskers in the box plots are long because no outliers
have been removed.  
We remark the ensemble average lies above the median because
it is dominated by a the few trajectories in 
${\cal E}$ that lead to errors at time $T$ with large magnitudes 
similar to Figure~\ref{epathsb}
while the geometric mean and median are much closer.
Since the median and geometric mean are more robust in the presence 
of extreme variations, the $\delta$ which minimizes the 
error shall be identified using these statistics.

We now address the difficulty that more frequent observations 
leads to lack of synchronization between the approximating
solution and reference solution by introducing a relaxation 
parameter into the
data assimilation algorithm given by Definition \ref{direct}.
To mitigate the difficulties that occur when inserting observations 
of the reference solution too frequently into the calculation,
it intuitively makes sense to insert those frequent observations 
less forcefully.  
We therefore
modify the discrete-in-time data assimilation
algorithm 
studied in \cite{COT2019}
by introducing a 
parameter $\kappa$ which will depend on $\delta$.
In details 
let $\kappa>0$ and 
replace the recurrence defining $u_{n+1}$ in Definition \ref{nokappa} by 
\begin{align}\label{relaxed}
	u_{n+1}&=
		(I-\kappa J) S(t_{n+1}-t_n)u_n+ \kappa J U(t_{n+1}).
\end{align}
Taking $\kappa=1$ recovers the original algorithm while 
$\kappa<1$ underrelaxes the system. 
Note that smaller values of $\kappa$ lead to 
proportionately smaller updates to the 
predicted state.

We remark that $\kappa$ plays a role similar to the term 
used by the Bayesian update for the ensemble Kalman filter 
which balances the tradeoff between 
the accuracy of the predicted state and the quality of the observations 
(see, for example, Evensen \cite{Evensen2007}
or Grudzien and Bocquet \cite{Grudzien2023}).
The context considered here, however, is that of noise-free 
observations with exactly known dynamics.
Thus, it is not noise 
and model error which lead us to introduce $\kappa$ but the
balance between 
lack of orthogonality in a rough interpolant
and smoothing by the dynamics.
Although noise and model error are important for applications,
the present research studies only the noise-free case to focus on 
how $\kappa$ and $\delta$ are related.

The role $\kappa$ plays in \eqref{relaxed} may be 
compared as follows to the role of $\mu$ in
the continuous nudging approach of \cite{AOT2014} given by	
\begin{equation}\label{cda}
{dv\over dt} +\nu Av +B(v,v)=f+\mu J(U-v),
\words{where}
v(t_0)=u_0.
\end{equation}
Note the effect of the nudging over the
interval $[t_n,t_{n+1}]$ is approximately
\begin{equation}\label{heursitica}
	\int_{t_n}^{t_{n+1}} \mu J(U-v) dt
		\approx \delta \mu J\big(U(t_n)-v(t_n)\big).
\end{equation}
Since at the beginning of the same interval
\eqref{relaxed} kicks the approximation by
$$
		-\kappa JS(t_n-t_{n-1})u_{n-1}+\kappa JU(t_n) =
\kappa J\big(U(t_n)-u(t_n^-)\big),
$$
it is natural to identify $\kappa$ with $\delta\mu$.
The numerical and theoretical analysis of this identification
is the main focus of the present paper.

After setting $\kappa=\mu\delta$ it follows that $\kappa\to 0$ as
$\delta\to 0$.  Thus, the more frequent the observations
the smaller the update based on them.
Our research begins by studying
how the optimal value of $\kappa$ depends
on $\delta$ numerically.  These simulations suggest 
$\kappa=\mu\delta$ is reasonable.
To further this connection between $\kappa$ and $\delta$ we then
prove analytically that the $\kappa$-relaxed discrete-in-time 
data assimilation algorithm 
with $\kappa=\delta\mu$
converges to the continuous nudging approach as $\delta\to 0$.

In particular, our main theoretical result is 
\begin{theorem}\label{mainresult}
Given a reference solution $U$ lying on the global
attractor of the incompressible two-dimensional
Navier--Stokes equations,
let $u$ be the approximating solution obtained
by~\eqref{relaxed} with $\kappa=\mu\delta$ and 
let $v$ be obtained by \eqref{cda}.
Then, for every $T>t_0$ it follows that
$$
	\sup\big\{\|u(t)-v(t)\|:t\in[t_0,T]\big\}\to 0
	\wwords{as}\delta \to 0.$$
Moreover, the above result holds for every $J=P_\lambda P_H I_h$ 
independently of
whether $h$ and $\lambda$ have been chosen so $\|u(t)-U(t)\|\to 0$ or even 
$\|v(t)-U(t)\|\to 0$ as $t\to \infty$.
\end{theorem}%
\noindent

The above result allows us to take $\delta\to 0$; however, we were
unable to obtain suitable choices of $h$, $\lambda$ and $\mu$ such
that $|u(t)-U(t)|\to 0$ for all $\delta$ sufficiently small.
Instead we have
\begin{theorem}\label{partialresult}
There exists $h$, $\lambda$ and $\mu$ such that
for every $\varepsilon>0$ there corresponds $T$ large enough and 
$\delta_0>0$ such that
$$
	\|u(T)-U(T)\|<\varepsilon \wwords{for all} 0<\delta<\delta_0
$$
and every reference solution $U$ lying on the global attractor.
\end{theorem}

The above result is a consequence 
of Theorem \ref{mainresult} combined with Theorem 2 in \cite{AOT2014} 
using the triangle inequality.  We state the proof here to
illustrate the limits of our present theory.
\begin{proof}[Proof of Theorem~\ref{partialresult}.]
From Theorem 2 in \cite{AOT2014} there exists $h$, $\lambda$ 
and $\mu$ such that
$$
	\|v(t)-U(t)\|\to 0\wwords{as} t\to\infty.
$$
Let $\varepsilon>0$.
Therefore, there exists $T$ large enough such that
$$
	\|v(t)-U(t)\|<\varepsilon/2\wwords{for all} t\ge T.
$$
Now choose $\delta_0>0$ in Theorem~\ref{mainresult} such that
$$
	\sup\big\{\|u(t)-v(t)\|:t\in[t_0,T]\big\}\le \varepsilon/2
	\wwords{for all} 0<\delta<\delta_0.
$$
It follows that
$$
	\|u(T)-U(T)\|\le |u(T)-v(T)|+|v(T)-U(T)|< \varepsilon.
$$
This finishes the proof.
\end{proof}

This paper is organized as follows.  
In Section~\ref{num} we study how the optimal value of $\kappa$ depends
on $\delta$ numerically and find that $\kappa=\mu\delta$ is 
reasonable.
Section~\ref{theo} begins with some theoretical preliminaries 
and ends with the proof of
Theorem~\ref{mainresult}.
Concluding remarks and directions for future work
appear in Section \ref{conclusion}.
Also included is 
Appendix~\ref{Tapriori} which provides explicit bounds on $v$
in terms of $T$ that are
used in Section~\ref{theo}.

\section{Numerical Results}\label{num}
This section describes the computations
which appear in the introduction.  We then present numerical
evidence that taking the relaxation parameter 
$\kappa$ proportional to $\delta$ results in numerical 
synchronization of the approximating solution with the reference solution
over a wide range of values for $\delta$ where $\delta\le 0.5$ is not too large.

All simulations are performed for a $2\pi$-periodic 
domain discretized on a $512\times 512$ spatial grid
using a dealiased spectral method with 115599 active
Fourier modes.  Integration in time was by means of
the fourth-order exponential time-differencing
method introduced by Cox and Matthews \cite{Cox2002}
with step size $\Delta t=0.0078125$.
To avoid loss of precision the coefficients for
the time-stepping method were obtained as in Kassam and Trefethen 
\cite{Kassam2005}.

The computer program that performs the data assimilation experiments
presented in this paper is written in the C 
programming language by the authors and compiled with GCC~\cite{GCC}.
This same code was previously used in~\cite{CO2022} to study the effects
of noise in time-delay nudging.
An overview of how this program works follows.
The MPICH MPI library~\cite{MPICH} distributes the 
computation
of the reference and corresponding approximating solutions
among the computational nodes.  Each MPI rank is
threaded using OpenMP and the Fourier transforms
performed using the FFTW library of Frigo and Johnson~\cite{Frigo2005}.
While running, our simulations compute the reference solution using 
MPI rank 0 and every $m$ time steps where $m\Delta t=\delta$
send the observational measurements of $U$ to the remaining ranks.
Those ranks then construct for different values of $\kappa$ 
the corresponding approximating solutions.

We remark the one-way flow 
of observations from rank 0 to the other MPI ranks mimic real-world 
data gathering.  This one-way flow of information also means 
network communication 
latency does not significantly affect the performance of our simulations.
As a result the data-assimilation experiments could be
parallelized across a loosely-coupled cluster 
of available machines
using relatively low-speed gigabit Ethernet.

We now describe further details of the data assimilation experiments.
The observational measurements are given by local spatial averages
around the points 
$$
	p_i=\big(
		\lfloor x_{i,1}(256/\pi)+0.5\rfloor,
		\lfloor x_{i,2}(256/\pi)+0.5\rfloor
		\big)
$$
on a $512\times 512$ grid where
$$
x_i\in \big\{\, \big((2n_1+1)\pi/9,(2n_2+1)\pi/9\big) 
	: n_1,n_2=0,\ldots,8\,\big\}.
$$
The spatial averages are computed as 
$$
	U_i\approx {1\over |{\cal B}_i|}\sum_{ p\in{\cal B}_i}
	 U(p\pi/256)
\wwords{where}
		{\cal B}_i=\big\{\,p\in\Z^2: |p-p_i|^2\le 24\,\big\}.
$$
The condition $|p-p_i|^2\le 24$ corresponds to an averaging radius
of about $\sqrt{24}\pi/256\approx 0.06$ in the $2\pi$-periodic domain.
Note there are 81 local averages and $|{\cal B}_i|=69$ grid points in 
each average.  Since the grid consists of $512^2=262144$ points 
total, the local averages cover about $2.13$ percent of the physical domain.
This percentage may be interpreted as 
how well
the spatial 
averages reflect point observations of the velocity fields.

As the MPI ranks used for computing an approximating solution receive
the observational data it is interpolated as
$$
	JU(p)=P_\lambda P_H\sum_{i=1}^{81} U_i \chi_i(p)
$$
where $\chi_i$ is the characteristic function for the square centered 
at $p_i$.  

To identify the relationship between $\kappa$ and $\delta$ we performed
data-assimilation experiments for sampling intervals of 
\begin{equation}\label{dofb}
\delta={m\Delta t} \wwords{where} 
		m=\Big\lfloor 228 \Big({3\over 4}\Big)^p\Big\rfloor
\wwords{for} p=0,1,\ldots,17
\end{equation}
and the relaxation parameters
\begin{equation}\label{kofq}
	\kappa=\Big({3\over 4}\Big)^{q/3} 
\wwords{for which} \kappa\in \big[0.0056,5\delta\big]
\wwords{and} q\in\Z.
\end{equation}
In all we consider 873 choices for $\delta$ and $\kappa$.
For each choice multiple simulations
were performed corresponding to 
an ensemble ${\cal E}$ of 500 different reference solutions $U$ lying on the global attractor.
This resulted in a total of $436500$ data assimilation experiments.

\begin{figure}[h!]
    \centerline{\begin{minipage}[b]{0.75\textwidth}
    \caption{\label{kappa}%
    The geometric mean of
    $\|U(T)-u(T)\|_{L^2}$
	at $T=t_0+1024$
    for varying values of $\delta$ and $\kappa$.}
    \end{minipage}}
    \centerline{
        \includegraphics[height=0.4\textwidth]{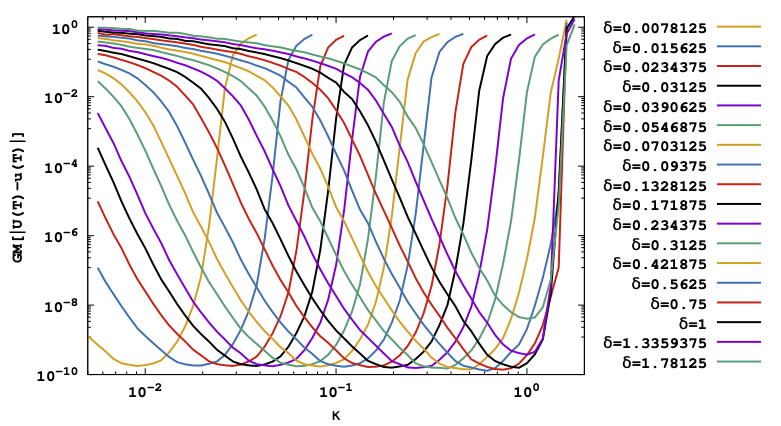}
    }
\end{figure}

The ensemble ${\cal E}$ was obtained
from the long-term evolutions of randomly-chosen points in
phase space.
In particular $U\in{\cal E}$ implies $U(t_0)=U_0$ where $U_0=S(10240) Z_0$
and $Z_0$ is a random velocity field
with energy spectrum similar to a point on the attractor.
Although $Z_0$ is unlikely to be on the global attractor, 
the evolution time of 10240 is long enough that large
scale structures have appeared in the velocity field and further
undergone more than 300 large eddy turnovers.  
We remark that the energetics of each flow appear to enter into
a statistically-steady state while continuing to undergo 
complex fluctuations.
This suggests the initial conditions $U_0$ for each of the
trajectories in ${\cal E}$ lie very near the global attractor and
that this attractor is nontrivial.

After every $3200$ time steps the reference and
approximating solutions were compared.  After 131072
time steps 
the final
value of $|U(T)-u(T)|$ at $T=t_0+1024$ is returned.
We note the available hardware---mostly Xeon 
and Epyc servers---typically 
performed from 10 to 30 time steps per second
depending on the exact processor and memory speed.
On average each run about took two hours.
To facilitate scheduling the runs were 
partitioned into batches that computed the reference solution
and corresponding approximating solutions 
for seven parameter choices at a time.

\begin{figure}[h!]
    \centerline{\begin{minipage}[b]{0.75\textwidth}
    \caption{\label{kdelta}%
    Box plots (no outliers removed) depicting the values 
	of $\kappa$ that minimize $\|U(T)-u(T)\|_{L^2}$
    when $T=t_0+1024$
    for each trajectory in the ensemble.
    Note $\kappa=\mu\delta$
    was fitted for $\delta\le 0.5$.
    }
    \end{minipage}}
    \centerline{
        \includegraphics[height=0.4\textwidth]{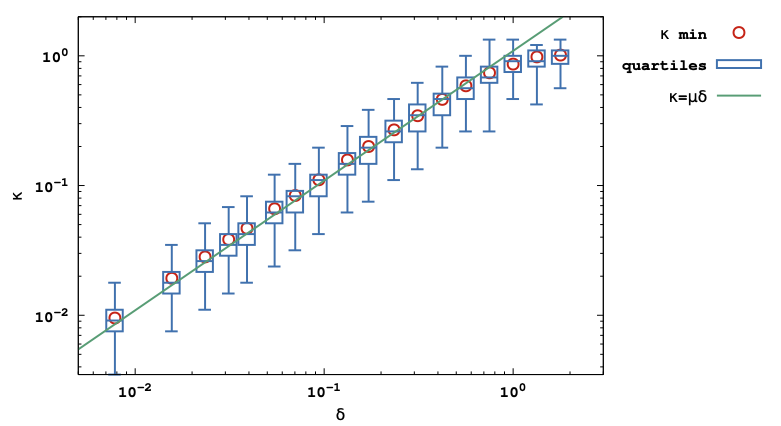}
    }
\end{figure}

Figure \ref{kappa} plots
the geometric mean of the error at time $T$ 
for each $\delta\in [0.0078125,1.78125]$ 
as a function of $\kappa$ where
$\kappa\in [0.0056,5\delta]$.
Note that additional runs when $\delta=0.0078125$ were performed
with $\kappa\in \big[0.00317,0.0056\big)$ to ensure that small 
enough values of $\kappa$ were considered 
for this smallest value of $\delta$.

The value $\kappa_{\rm min}$ of 
$\kappa$ that minimizes the geometric mean 
along the corresponding curve for each $\delta$ 
in Figure~\ref{kappa} 
was obtained using a least-squares quadratic fit of 
points near the minimum.  
To illustrate how $\kappa_{\rm min}$ depends on $\delta$ 
the points $(\delta,\kappa_{\rm min})$ are plotted 
in Figure~\ref{kdelta} as circles.
Note for $\delta\le 0.5$ the relationship
between $\kappa_{\rm min}$ and $\delta$ appears linear.
This linearity is our primary numerical result and the focus of 
the theory in the following section.

Figure~\ref{kdelta} further characterizes the relationship
between $\kappa$ and $\delta$ by considering each reference
trajectory $U\in{\cal E}$ separately to obtain 500 different
minimizing values of $\kappa$ for each $\delta$.
The distribution of the minimizing values of $\kappa$ are then displayed
in quartiles for each value of $\delta$ considered.
Finally a least squares fit of $\kappa=\mu\delta$ for $\delta\le 0.5$
was performed to obtain $\mu\approx 0.966$.
We remark $\mu$ is not universal but a tuning
parameter that depends on the viscosity $\nu$, body forcing $f$, 
spatial domain and the interpolant observable $J$ among others.
In particular, the value of $\mu$ obtained above is only 
coincidentally close to one.  What is important for the present
research is that
the relationship between $\kappa$ and $\delta$ appears linear.

\section{Theoretical Results}\label{theo}

In this section we prove our main theoretical result 
given as Theorem~\ref{mainresult} in the introduction.
Begin by translating the {\it a priori\/} bounds
of Theorem~\ref{rhoapriori} into bounds on the 
interpolated observations represented by $JU$.

\begin{proposition}\label{JUbound}
\noindent Suppose $U$ lies on the global attractor of \eqref{nse1}
and $f\in V$.  There exists constants $\rho_J$ and $\rho_K$ depending
on $\lambda$, $\nu$, $\|f\|$ and $h$ such that
\begin{equation}\label{Jrho}
	|J U| \le \rho_J\wwords{and}
	\|J U\| \le \rho_K.
\end{equation}
\end{proposition}
\begin{proof}  Estimate using \eqref{Jinterp} followed
by \eqref{Urhov} as
$$
	|JU|=|U-JU|+|U|
	\le (\lambda^{-1}+c_1h^2)^{1/2} \|U\| + |U|
$$
to obtain $\rho_J= (\lambda^{-1}+c_1h^2)^{1/2}\rho_1+\rho_0$.
Similarly use \eqref{vcont} to estimate
$$
	\|J U\|=\|U-JU\|+\|U\|
	\le (1+\lambda c_1h^2)^{1/2} \|U\| + \|U\|.
$$
Thus, $\rho_K=(1+\lambda c_1h^2)^{1/2}\rho_1+\rho_1$.
\end{proof}

Our estimates make use of the following version of Theorem 5 from 
\cite{AOT2014} that removes the condition on $h$ and $\mu$ 
but applies only on a finite time interval $[t_0,T]$.
\begin{theorem}\label{thmvbound}
Suppose $f\in V$ and $v_0\in {\cal D}(A)$.  Let $v$ be the solution obtained
from \eqref{cda}.  There exists
bounds $R_\alpha$ and $\widetilde R_\alpha$ depending 
on $T$, $t_0$, $\nu$, $\|f\|$ and $v_0$ such that
\beq\label{vbound1}
	\|v\|_\alpha\le R_\alpha
	\wwords{for} t\in[t_0,T]
	\wwords{with} \alpha=0,1,2
\eeq
and 
\beq\label{vbound2}
	\Big(\int_{t_n}^{t_{n+1}}\|v\|_\alpha^2dt\Big)^{1/2}
		\le \widetilde R_\alpha
	\words{for} [t_n,t_{n+1}]\subseteq[t_0,T]
	\words{with} \alpha=1,2,3.
\eeq
\end{theorem}%
\noindent 
The proof of Theorem \ref{thmvbound}
is provided in Appendix \ref{Tapriori}.
We note here that the smoothness of the interpolant $J$ as well as $f$
constrains the regularity of $v$.

We remark 
that $v_0=u_0$
in \eqref{cda}
where $u_0=JU(t_0)$.  Since $J=P_\lambda P_H I_h$ then
$P_\lambda v_0=v_0$.
Consequently
\eqref{psmooth} followed by Proposition~\ref{JUbound} yields the bound
$$|Av_0|=\|P_\lambda v_0\|_2\le\lambda |v_0|
	=\lambda|JU(t_0)|\le\lambda\rho_J.$$
In particular, we may assume $v_0\in{\cal D}(A)$
in the proof of Theorem~\ref{mainresult}.
Moreover, since the reference solution $U$ lies on
the global attractor $\A$ we may further assume the bounds 
$R_\alpha$ and $\widetilde R_\alpha$
to be independent of $U_0$.

Now use \eqref{Jinterp} 
and \eqref{vcont} 
as in the proof of Proposition~\ref{JUbound}
to obtain bounds on $|Jv|$ and $\|Jv\|$.
Thus, there exists $R_J$ and $R_K$ further
depending on $h$ and $\lambda$ such that
\begin{equation}\label{JR}
	|J v|\le R_J\wwords{and}\|J v\|\le R_K.
\end{equation}

A number of inequalities involving the non-linear term
are summarized in Chapter II Appendix A of 
Foias, Manley, Rosa and Temam \cite{Foias2001}.
We recall here for reference
\begin{equation}\label{bbthree}
|(B(u,v),w)|\le c|u|^{1/2}\|u\|^{1/2}\|v\|^{1/2}|Av|^{1/2}|w|
\end{equation}
for $u\in V$, $v\in {\cal D}(A)$ and $w\in H$
as well as
\begin{equation}\label{bbfive}
|(B(u,v),w)|\le c|u|^{1/2}|Au|^{1/2}\|v\||w|
\end{equation}
for $u\in {\cal D}(A)$, $v\in V$ and $w\in H$.
Note that the above inequalities hold for a suitable constant $c$ 
depending on the domain $\Omega$ which in our case is $L$-periodic.

It follows from \eqref{bbthree} that
\begin{equation}\label{bbone}
|B(u,v)|
       \le c |u|^{1/2} \|u\|^{1/2} \|v\|^{1/2} |Av|^{1/2}
	\wwords{for} u\in V \hbox{ and } v\in {\cal D}(A)
\end{equation}
and from \eqref{bbfive} that
\begin{equation}\label{bbtwo}
|B(u,v)|
   \le c |u|^{1/2} |Au|^{1/2} \|v\|
\wwords{for} u\in {\cal D}(A) \hbox{ and } v\in V.
\end{equation}
Again from \eqref{bbone} it follows that
\begin{equation}\label{buvnorm}
	\|B(u,v)\|\le c\|u\|^{1/2}|Au|^{1/2}\|v\|^{1/2}|Av|^{1/2}
		+c|u|^{1/2}\|u\|^{1/2}|Av|^{1/2}|A^{3/2}v|^{1/2}
\end{equation}
for $u\in{\cal D}(A)$ and $v\in{\cal D}(A^{3/2})$.
Setting $u=v$ and interpolating then yields 
\begin{equation}\label{bvnorm}
	\|B(u,u)\|\le
		c|u|^{1/2}\|u\|^{1/2}|Au|^{1/2}|A^{3/2}u|^{1/2}
\wwords{for} u\in{\cal D}(A^{3/2}).
\end{equation}
For simplicity the $c$ appearing above and in \eqref{logestB} below
will denote a single constant chosen suitably large so all of the 
inequalities hold.

We shall also employ an estimate appearing in 
\cite{CO2022},
Foias, Mondaini and Titi~\cite{foias2016discrete}
and Titi~\cite{titi1987criterion}
based on the Brezis--Gallou\"et inequality
which we state here as
\begin{proposition}
If $v$ and $w$ are in ${\cal D}(A)$ then 
\begin{equation}\label{logestB}
	\big|\big(B(w,v)+B(v,w),Aw\big)\big| 
	\le c \|w\| \|v\|\bigg( 
	1+\log{|Av|\over \lambda_1^{1/2}\|v\|}\bigg)^{1/2} |Aw|,
\end{equation}
where $c$ is a non-dimensional constant depending only on the domain.
\end{proposition}

Fix $T>t_0$.  
To study the limit 
on the interval $[t_0,T]$ 
as the time between observations $\delta\to 0$
define $t_n=t_0+\delta n$ for
$n=0,\ldots,N$ where
$N$ is the greatest integer such that $t_N\le T$.
Assume $\delta\le\dmax$ where $\dmax\le T-t_0$; otherwise,
if $\delta>T-t_0$, then no observational measurements would occur 
during the interval $[t_0,T]$ and there would be no data to assimilate.
Upon substituting $\kappa=\delta\mu$ 
into~\eqref{relaxed} and using the fact that $t_{n+1}-t_n=\delta$, the
relaxed version of Definition~\ref{direct} becomes
\begin{definition}\label{relaxudef}
Assume $\delta\le T-t_0$ and define $t_n=t_0+\delta n$ for
$n=0,\ldots,N$ where
$N$ is the greatest integer such that $t_N\le T$.
The relaxed direct-insertion method is defined as the
approximating solution $u(t)$ for $t\in [t_0,T]$ 
given by
$$
    u(t)=
\begin{cases}
	S(t-t_n)u_n &\hbox{for}\quad t\in[t_n,t_{n+1}) 
		\words{with} n=0,\ldots,N-1,\\
	S(t-t_N)u_N &\hbox{for}\quad t\in[t_N,T]
\end{cases}
$$
where
$$
\left\{\begin{aligned}
    u_0&=J U(t_0),\nonumber\\
	u_{n+1}&=S(\delta)u_n +\delta\mu J
		\big\{U(t_{n+1})-S(\delta)u_n\big\}.
\end{aligned}\right. 
$$
\end{definition}

In preparation to prove Theorem~\ref{mainresult}
let $w_n=v(t_n)-u_n$ and $w(t)=v(t)-u(t)$.
Here $v(t)$ is the approximating solution
obtained by \eqref{cda}.
Definition \ref{relaxudef} implies
$$
	u(t_{n+1}^-)=
	 	\lim_{t\nearrow t_{n+1}} u(t)=
		S(\delta)u_n.
$$
Therefore
\begin{align*}
w_{n+1}
	&=v(t_{n+1})-(1-\delta\mu J)u(t_{n+1}^-)-\delta\mu JU(t_{n+1})\\
	&=w(t_{n+1}^-)-\delta\mu J(U(t_{n+1})-v(t_{n+1})+w(t_{n+1}^-))\\
	&=(I-\delta\mu J)w(t_{n+1}^-)-\delta\mu J(U(t_{n+1})-v(t_{n+1})).
\end{align*}
Consequently $w$ satisfies
\beq\label{cda1}
\frac{dw}{dt}+\nu Aw+B(v,w)+B(w,v)-B(w,w)=\mu J(U-v)
	\words{on}  t\in[t_n,t_{n+1})
\eeq
with initial condition
\beq\label{wn1}
	w_{n}=(I-\delta\mu J)w(t_{n}^-)
	-\delta\mu J(U(t_{n})-v(t_{n}))\wwords{for}n>0
\eeq
and $w_0=0$ since we use the same initial state for both
data assimilation methods.  For notational convenience 
we define $w(t_0^-)=0$ so that \eqref{wn1} also
holds when $n=0$.

Continue with the following estimate on the integral of $L^2$ norm of $Aw$.

\begin{theorem}\label{lemmaAw}
There exists $M_0$, $M_1$ and $\widetilde M_2$ depending only on 
$T$, $t_0$, $\nu$, $\|f\|$ and the domain $\Omega$ such that
$$
	|w|\le M_0\wwords{and}
	\|w\|\le M_1
\wwords{for all} t\in [t_0,T]
$$
and
\begin{equation}\label{tildeM2bound}
\Big(\int_{t_n}^{t_{n+1}} | Aw|^2dt\Big)^{1/2}\le \widetilde M_2 
	\wwords{for} 
[t_n,t_{n+1}]\subseteq [t_0,T].
\end{equation}
\end{theorem}
\begin{proof}
Let $N$ be the greatest integer such that $t_N\le T$.
By Theorem \ref{thmvbound} 
there exists $R_\alpha$ 
such that
$\|v\|_\alpha\le R_\alpha$
for $t\in[t_0,T]$ and $\alpha=0,1,2$.
Since the reference trajectory $U$ 
lies on the global attractor and $v_0=JU(t_0)$ than these bounds 
can be taken independent of $v_0$.
Multiply \eqref{cda1} by $Aw$ and use the fact 
that $(B(w,w),Aw)=0$ in two-dimensional periodic domains 
to obtain
\beq\label{cda12}
\frac{1}{2}\frac{d\|w\|^2}{dt}
	+\nu |Aw|^2+(B(v,w),Aw)+(B(w,v),Aw)=\mu (J(U-v),Aw).
\eeq
From \eqref{logestB} it follows that
\begin{align}
|(B(v,w),Aw)&+(B(w,v),Aw)|\le 
	c\|w\|\|v\|\bigg(
		1+\log\frac{|Av|}{\lambda_1^{1/2}\|v\|}\bigg)^{1/2}|Aw|
		\nonumber\\
	\label{B2est}
	&\le cR_1^{1/2}R_L^{1/2}\|w\||Aw|
	\le \frac{c^2R_1R_L}{\nu}\|w\|^2+\frac{\nu}{4}|Aw|^2,
\end{align}
where
$$
	R_L=R_1\bigg(1+\log {R_2\over \lambda_1^{1/2}R_1}\bigg).
$$
Also
\beq\label{m1est}
	\mu (J(U-v),Aw)\le \frac{\mu^2}{\nu} |J(U-v)|^2+\frac{\nu}{4}|Aw|^2
	\le \frac{\mu^2}{\nu} \ca^2+\frac{\nu}{4}|Aw|^2
\eeq
where
$\ca=\rho_J+R_J$.

Combining \eqref{B2est} and \eqref{m1est} with \eqref{cda12} yields 
\beq\label{cda123}
	\frac{d\|w\|^2}{dt}+\nu |Aw|^2
	\le\frac{2c^2R_1R_L}{\nu}\|w\|^2+\frac{2\mu^2}{\nu} 
		\ca^2\le \cb\|w\|^2+\cc.
\eeq
Here $\cb=2c^2R_1R_L/\nu$ and 
$\cc=2\mu^2 \ca^2/\nu$.
Now multiply \eqref{cda123} by $e^{-\cb t}$ and integrate
from $t_n$ to $t$ where $t\in [t_n,t_{n+1})$ for $n=0,1,\ldots,N-1$
and $t\in[t_N,T]$ for $n=N$.  
Thus,
\begin{align}\label{integrated}
	\|w(t)\|^2&+
		\nu \int_{t_n}^t e^{\cb(t-s)} |Aw(s)|^2 ds
	\le e^{\cb(t-t_n)}\|w_n\|^2+{\cc\over \cb}\big(e^{\cb(t-t_n)}-1\big).
\end{align}
We obtain after applying
\eqref{vcont} that
\begin{align*}
	\|w_n\|
	&=\|(I-\delta\mu J)w(t_n^-)
	    -\delta\mu J(U(t_n)-v(t_n))\|\cr
	&\le  |1-\delta\mu|\|w(t_n^-)\|
		+\delta\mu\|(I-J)w(t_n^-)\|
	    +\delta\mu \|J(U(t_n)-v(t_n))\|.
\end{align*}
Therefore
\begin{equation}\label{wnest}
	\|w_n\|
	\le
		(1+\delta\mu\cd)\|w(t_n^-)\|
	    +\delta\mu\ce
\end{equation}
where $\cd=1+\sqrt{1+\lambda c_1h^2}$
and $\ce=\rho_K+R_K$.
Consequently, under the assumption $\delta\le \dmax$ 
it follows that
\begin{align*}
	\|w_n\|^2
	&\le (1+\delta\mu\cd)^2\|w(t_n^-)\|^2
		+2\delta\mu(1+\delta\mu\cd)\ce\|w(t_n^-)\|
		+\delta^2\mu^2\ce^2\cr
	&\le (1+\delta\mu\cd)^2\|w(t_n^-)\|^2
		+\delta\mu\big\{(1+\delta\mu\cd)^2 \|w(t_n^-)\|^2+\ce^2\big\}
		+\delta^2\mu^2\ce^2\cr
	&= (1+\delta\mu)(1+\delta\mu\cd)^2\|w(t_n^-)\|^2
		+\delta\mu(1+\delta\mu)\ce^2\cr
	&\le (1+\delta\mu\cf)\|w(t_n^-)\|^2 +\delta\mu\cg
\end{align*}
where 
$\cf=2\cd +\dmax\mu\cd^2 +(1+\dmax\mu\cd)^2$ 
and $\cg=(1+\dmax\mu)\ce^2$.

Substituting into \eqref{integrated} and dropping the integral
on the left results in
\begin{align}\label{tbound}
	\|w(t)\|^2\le e^{\cb(t-t_n)}
\big\{
	 (1+\delta\mu\cf)\|w(t_n^-)\|^2 +\delta\mu\cg
\big\}
		+{\cc\over \cb}\big(e^{\cb(t-t_n)}-1\big).
\end{align}
Subsequently for $n=0,1,\ldots,N-1$ taking the limit as $t\to t_{n+1}$ obtains
\begin{align*}\label{istep}
	\|w(t_{n+1}^-)\|^2
	&\le
	e^{\delta \cb}(1+\delta\mu\cf)
\|w(t_n^-)\|^2 
	+\delta\mu e^{\delta \cb}\cg+{\cc\over\cb}\big(e^{\delta\cb}-1\big).
\end{align*}
Setting $y_n=\|w(t_n^-)\|^2$ gives the inequality 
$y_{n+1}\le \alpha y_n+\beta$
where $\alpha=e^{\delta\cb}(1+\delta\mu\cf)$ and 
$\beta=\delta\mu e^{\delta\cb}\cg+(\cc/\cb)(e^{\delta\cb}-1)$.
By induction
\begin{align*}
	y_1&\le \alpha y_0+\beta\cr
	y_2&\le \alpha y_1+\beta \le \alpha^2 y_0+(1+\alpha)\beta\cr
	&\,\,\,\vdots\cr
	y_n&\le \alpha^n y_0+\textstyle\sum_{j=0}^{n-1}\alpha^j\beta
		\le \alpha^n y_0+{\alpha^n-1\over\alpha-1}\beta.
\end{align*}

Using the fact that
$x\le e^x-1\le x e^x$ for $x\ge 0$ we obtain
$$
	\alpha-1=e^{\delta\cb}-1+\delta\mu e^{\delta\cb}\cf
		\ge \delta (\cb+\mu e^{\delta\cb}\cf),
$$
$$
	\alpha^n =e^{n\delta\cb}(1+\delta\mu\cf)^n
		\le e^{(t_n-t_0)(\cb+\mu\cf)}
$$
and
$$
	\beta=\delta\mu e^{\delta\cb}\cg+(\cc/\cb)(e^{\delta\cb}-1)
		\le \delta e^{\delta\cb}(\mu\cg+\cc).
$$
Consequently, for $n=1,2,\ldots,N$ it follows from $\delta\le\dmax$ that
\begin{align*}
	\|w(t_n^-)\|^2
	&\le e^{(t_n-t_0)(\cb+\mu\cf)} \|w(t_0^-)\|^2
		+{e^{(t_n-t_0)(\cb+\mu\cf)}-1
		\over \cb+\mu e^{\delta\cb}\cf}
		e^{\delta\cb}(\mu\cg+\cc)\cr
	&\le e^{2\dmax(\cb+\mu\cf)} \|w(t_0^-)\|^2
		+{e^{2\dmax(\cb+\mu\cf)}-1
		\over (\cb+\mu\cf)}
		e^{\dmax\cb}(\mu\cg+\cc).
\end{align*}

Since $\|w(t_0^-)\|=0$, there is a constant $\ch$ independent
of $\delta$ such that 
$$
	\|w(t_n^-)\|\le\ch\wwords{for} n=1,\ldots,N.
$$
It follows from \eqref{wnest} that
$$
	\|w_n\|
	\le(1+\delta\mu\cd)\|w(t_n^-)\| +\delta\mu\ce
	\le (1+\dmax\mu\cd)\ch+\dmax\mu\ce
$$
and so there is also a constant $\ci$ independent of $\delta$ such that
$$
	\|w_n\|\le\ci\wwords{for} n=1,\ldots,N.
$$

Now use \eqref{integrated} to show there is $M_1$ such that 
$\|w(t)\|\le M_1$ for $t\in [t_0,T]$.
Note that
$$
	M_1^2 = 
	e^{\cb \dmax}\ci^2+{\cc\over \cb}\big(e^{\cb\dmax}-1\big).
$$
The Poincar\'e inequality 
\eqref{poincare} stated as
$\lambda_1 |w|^2 \le \|w\|^2$ immediately implies
$$
	|w(t)|\le \lambda_1^{-1/2}\|w(t)\|\le M_0
	\wwords{for} t\in[t_0,T]
$$
where $M_0=\lambda_1^{-1/2}M_1$.

To obtain \eqref{tildeM2bound} suppose $n=0,1,\ldots,N-1$ and 
integrate \eqref{cda123} from $t_n$ to $t$.
Then taking the limit as $t\to t_{n+1}$ yields
$$
	\|w(t_{n+1}^-)\|^2-\|w_n\|^2
	+\nu \int_{t_n}^{t_{n+1}}|Aw|^2dt
	\le C_2\int_{t_n}^{t_{n+1}}\|w\|^2dt+\delta C_3.
$$
Thus
$$
\nu \int_{t_n}^{t_{n+1}}|Aw|^2dt\le \|w_n\|^2+ 
	\dmax(C_2M_1^2+C_3)
$$
or equivalently 
$$
 \Big(\int_{t_n}^{t_{n+1}}|Aw|^2dt\Big)^{1/2}\le \widetilde M_2
\wwords{where} \widetilde M_2^2=
	{\ci^2\over \nu}+ \dmax{C_2M_1^2 +C_3\over\nu}.
$$
Finally, noting that $M_0$, $M_1$ and $\widetilde M_2$ depend only 
on $T$, $t_0$, $\nu$, $\|f\|$
and $\Omega$ finishes the proof.
\end{proof}
\begin{corollary}\label{Jcor} There exists constants $M_J$ and $M_K$
such that
$$
	|Jw|\le M_J\wwords{and} \|J w\|\le M_K
	\wwords{for all} t\in[t_0,T].
$$
\end{corollary}
\begin{proof}
This follows at once from Theorem~\ref{lemmaAw} applied to
\eqref{Jinterp} and \eqref{vcont} exactly analogous to the way 
\eqref{Jrho} and \eqref{JR} were obtained.
\end{proof}

The next step is to bootstrap Theorem~\ref{lemmaAw} to obtain
bounds which tend to zero as $\delta\to 0$.
In order to do this we first prove the following two lemmas.
\begin{lemma}\label{intdwbound}
There exists $M_D$ independent of $\delta$ such that
$$
	\int_{t_n}^{t_{n+1}} |w(t)-w_n| dt \le \delta^{3/2} M_D
	\wwords{for} [t_n,t_{n+1}]\subseteq [t_0,T].
$$
\end{lemma}
\begin{proof}
Estimate as
\begin{align*}
	\int_{t_n}^{t_{n+1}} &|w(t)-w_n|dt
	\le\int_{t_n}^{t_{n+1}}\int_{t_n}^{s} 
		\Big|\frac{dw(t)}{dt}\Big|dt\,ds
	\le\delta\int_{t_n}^{t_{n+1}}
		\Big|\frac{dw(t)}{dt}\Big|dt\\
	&\le \delta\int_{t_n}^{t_{n+1}}
		\big\{\nu|Aw|+|B(v,w)|+|B(w,v)|+|B(w,w)|
		+\mu|J(U-v)|\big\}.
\end{align*}
By the Theorem~\ref{lemmaAw} we have 
\beq\label{cda16a}
\begin{aligned}
	\int_{t_n}^{t_{n+1}}\nu|Aw|
	\le\nu\delta^{1/2}\Big(\int_{t_n}^{t_{n+1}}|Aw|^2\Big)^{1/2}
	\le\delta^{1/2}\nu{\widetilde M_2}.
\end{aligned}
\eeq
Using \eqref{bbone} we have
\begin{align*}
	\int_{t_n}^{t_{n+1}} |B(v,w)| &\le
	c\int_{t_n}^{t_{n+1}} |v|^{1/2}\|v\|^{1/2}\|w\|^{1/2}|Aw|^{1/2}
	\le c R_0^{1/2}R_1^{1/2}M_1^{1/2}
	\int_{t_n}^{t_{n+1}} |Aw|^{1/2}\cr
	&\le c R_0^{1/2}R_1^{1/2}M_1^{1/2} \delta^{3/4}
		\Big(\int_{t_n}^{t_{n+1}} |Aw|^2\Big)^{1/4}
	\le \delta^{3/4}c R_0^{1/2}R_1^{1/2}M_1^{1/2}\widetilde M_2^{1/2}
\end{align*}
and from \eqref{bbtwo} it follows
\begin{align*}
	\int_{t_n}^{t_{n+1}} |B(w,v)| &\le
	c\int_{t_n}^{t_{n+1}} |w|^{1/2}|Aw|^{1/2}\|v\|
	\le \delta^{3/4} cM_0^{1/2} R_1 \widetilde M_2^{1/2}.
\end{align*}
as well as
\begin{align*}
	\int_{t_n}^{t_{n+1}} |B(w,w)| &\le
	c\int_{t_n}^{t_{n+1}} |w|^{1/2}|Aw|^{1/2}\|w\|
	\le \delta^{3/4} cM_0^{1/2} M_1 \widetilde M_2^{1/2}.
\end{align*}
Finally,
$$
	\int_{t_n}^{t_{n+1}} \mu|J(U-v)|
		\le \delta\mu(\rho_J+R_J).
$$

Upon collecting the above estimates together we obtain
$$
	\int_{t_n}^{t_{n+1}} |w(t)-w_n|dt
		\le \delta^{3/2} M_D
$$
where
$$
	M_D=\nu\widetilde M_2
	+\dmax^{1/4}c \big\{
		R_0^{1/2}R_1^{1/2}M_1^{1/2}
		+M_0^{1/2} R_1+M_0^{1/2} M_1\big\}\widetilde M_2^{1/2}
		+\dmax^{1/2} \mu(\rho_J+R_J).
$$
Noting that $M_D$ is independent of $\delta$ finishes the proof.
\end{proof}
\begin{corollary}\label{intdwboundcor}
There exists $M_E$ independent of $\delta$ such that
$$  
    \int_{t_n}^{t_{n+1}} |w(t)-w(t_n^-)| dt \le \delta^{3/2} M_E
    \wwords{for} [t_n,t_{n+1}]\subseteq [t_0,T].
$$  
\end{corollary}
\begin{proof}
Since
$ w_n - w(t_n^-)
=-\delta\mu J\big(w(t_n^-)+U(t_n)-v(t_n)\big)
$
then
$$
    \int_{t_n}^{t_{n+1}} |w_n-w(t_n^-)| dt
	\le \delta^2\mu (M_J+\rho_J+R_J).
$$
It follows from Lemma \ref{intdwbound} that
$$
	\int_{t_n}^{t_{n+1}} |w(t)-w(t_n^-)| dt
	\le \delta^{3/2} M_D+\delta^2\mu (M_J+\rho_J+R_J)
	\le \delta^{3/2}M_E
$$
where $M_E= M_D+\dmax^{1/2}\mu (M_J+\rho_J+R_J)$.
\end{proof}

\begin{lemma}\label{intdvbound}
	There exists $\rho_D$ and $R_D$ independent of 
$\delta$ such that for $[t_n,t_{n+1}]\subseteq[t_0,T]$ holds
$$
\int_{t_n}^{t_{n+1}}\|J(U(t)-U(t_n))\|dt \le \delta^{3/2} \rho_D
\words{and}
\int_{t_n}^{t_{n+1}}\|J(v(t)-v(t_n))\|dt \le \delta^{3/2} R_D.
$$
\end{lemma}
\begin{proof}
Estimate using \eqref{Jsmooth}
\begin{align*}
	\int_{t_n}^{t_{n+1}}\|J(U(t)-U(t_n))\|dt
	&\le c_3 \int_{t_n}^{t_{n+1}}\|U(t)-U(t_n)\|dt
	\le c_3 \int_{t_n}^{t_{n+1}}\int_{t_n}^{s}
		\Big\|\frac{dU(t)}{dt}\Big\|dt\,ds\\
	&\le c_3 \delta \int_{t_n}^{t_{n+1}}
		\big\{\nu\|AU\|+\|B(U,U)\|+\|f\|\big\}
\end{align*}
By the Cauchy--Schwartz inequality and 
\eqref{Uintbound}
$$
	\nu \int_{t_n}^{t_{n+1}} \|AU\|
	\le \delta^{1/2} \nu
	\Big(\int_{t_n}^{t_{n+1}} \|AU\|^2\Big)^{1/2}\le \delta^{1/2}
		\nu\widetilde\rho_3.
$$
Similarly \eqref{bvnorm} implies
\begin{align*}
	\int_{t_n}^{t_{n+1}} \|B(U,U)\|
	&\le c 
	\int_{t_n}^{t_{n+1}} |U|^{1/2}\|U\|^{1/2}|AU|^{1/2}|A^{3/2}U|^{1/2}\cr
	&\le c 
	\int_{t_n}^{t_{n+1}} \rho_0^{1/2}\rho_1^{1/2}	
		|AU|^{1/2}|A^{3/2}U|^{1/2}
	\le \delta^{1/2} c
		(\rho_0\rho_1 \widetilde\rho_2\widetilde\rho_3)^{1/2}.
\end{align*}
Therefore,
$$
	\rho_D=c_3\big(\nu\widetilde\rho_3+
        c (\rho_0\rho_1
        \widetilde\rho_2\widetilde\rho_3)^{1/2}
        +\dmax^{1/2} \|f\|\big).
$$

The proof of the second inequality is similar except there is one 
more term which we estimate using \eqref{Jrho} and \eqref{JR} as
\begin{align*}
	\int_{t_n}^{t_{n+1}}&\|J(v(t)-v(t_n))\|dt
	\le c_3\int_{t_n}^{t_{n+1}}\|v(t)-v(t_n)\|dt
	\le c_3\int_{t_n}^{t_{n+1}}\int_{t_n}^{s}
			\Big\|\frac{dv(t)}{dt}\Big\|dt\,ds \cr
&\le c_3\delta\int_{t_n}^{t_{n+1}}
	\big\{\nu\|Av\|+\|B(v,v)\|+\|f\|+\mu\|J(U-v)\|\big\}
	\le \delta^{3/2} R_D
\end{align*}
where
$
	R_D=c_3\big(\nu\widetilde R_3+
        c (R_0R_1 \widetilde R_2\widetilde R_3)^{1/2}
        +\dmax^{1/2} (\|f\|
		+\mu\rho_K+\mu R_K)\big)
$.
\end{proof}

We are now ready to prove our main theoretical result
stated as Theorem \ref{mainresult} in the introduction.
Namely, under the above hypotheses we prove
\begin{theorem}
For every $T>t_0$ that
$$
    \sup\big\{\|w(t)\|:t\in[t_0,T]\big\}\to 0
    \wwords{as}\delta \to 0.
$$
\end{theorem}
\begin{proof}
Fix $T>t_0$. Let $t_n=t_0+\delta n$ where $\delta\le\dmax$
and $\dmax\le T-t_0$.
Let $N$ be the greatest integer such that $t_N\le T$
and consider the interval $[t_0,T]$.
Taking the inner product of \eqref{cda1} with $Aw$ as in the proof of
Theorem \ref{lemmaAw}
we obtain
\begin{equation}\label{zcda3pre1}
	\frac{1}{2}\frac{d\|w\|^2}{dt}
		+{1\over 2}\nu |Aw|^2\le {c^2 R_1R_L\over 2\nu}\|w\|^2
		+\mu (J(U-v),Aw).
\end{equation}
Note the factor-of-two improvement in the coefficient
of $\|w\|^2$ on the right comes from a different 
application of Young's inequality.
In particular, rather than \eqref{B2est}
we estimate
$$
|(B(v,w),Aw)+(B(w,v),Aw)|
    \le \frac{c^2R_1R_L}{2\nu}\|w\|^2+\frac{\nu}{2}|Aw|^2.
$$
Upon applying the Poincar\'e inequality \eqref{poincare}
written as $\lambda_1\|w\|^2\le |Aw|^2$ 
we obtain
\begin{equation}\label{zcda3}
	\frac{d\|w\|^2}{dt}
		\le \cl\|w\|^2
		+2\mu (w,AJ(U-v)).
\end{equation}
where $\cl=c^2R_1R_L/\nu-\nu\lambda_1$.

Recall that $w(t_n)=w_n$.
The definition of $w_n$ given in \eqref{wn1} implies
\begin{align*}
	w_n
	&=(1-\delta\mu)w(t_n^-)+\delta\mu(I-J)w(t_n^-)-\delta\mu J(U(t_n)-v(t_n)).
\end{align*}
Therefore,
\begin{align*}
\|w_n\|^2
	&=(1-\delta\mu)^2\|w(t_n^-)\|^2+\delta^2\mu^2\|(I-J)w(t_n^-)
		-J(U(t_n)-v(t_n))\|^2\\
		&\qquad +2\delta\mu\big(A^{1/2}
			(1-\delta\mu)w(t_{n}^-),A^{1/2}
			(I-J)w(t_n^-)
			\big)\cr
		&\qquad -2\delta\mu(1-\delta\mu)\big(
			A^{1/2}w(t_{n}^-),
			A^{1/2}J(U(t_{n})-v(t_{n}))
			\big)\cr
	&\le (1+\delta\mu\cj)\|w(t_n^-)\|^2 
		-2\delta\mu\big(
			w(t_{n}^-),
			AJ(U(t_{n})-v(t_{n}))
			\big)
		+\delta^2\mu^2\ck
\end{align*}
where
$$
	\cj=2+\dmax\mu+2c_3(1+\dmax\mu)
$$
and
$$
	\ck=(M_1+M_K+\rho_K+R_K)^2+2 M_1(\rho_K+R_K).
$$
Note for the above estimates we have used
\begin{equation*}
	(1-\delta\mu)^2\le 1+\delta\mu(2+\dmax \mu),
\end{equation*}
\begin{equation*}
	\delta^2\mu^2\|(I-J)w(t_n^-)
        -J(U(t_n)-v(t_n))\|^2
	\le \delta^2\mu^2(M_1+M_K+\rho_K+R_K)^2,
\end{equation*}
\begin{align*}
	2\big(A^{1/2}&(1-\delta\mu)w(t_{n}^-),
		A^{1/2}(I-J)w(t_n^-)\big)
	\le 
		2(1+\dmax\mu)\|w(t_n^-)\|\|(I-J)w(t_n^-)\|\cr
	&\le
		2(1+\dmax\mu) (1+\lambda c_1h^2)^{1/2} \|w(t_n^-)\|^2
	\le 
		2c_3(1+\dmax\mu)\|w(t_n^-)\|^2
\end{align*}
and
\begin{align*}
	2\delta^2\mu^2&\big(
			A^{1/2}w(t_{n}^-),
			A^{1/2}J(U(t_{n})-v(t_{n}))
			\big)\le 2\delta^2\mu^2 M_1(\rho_K+R_K).
\end{align*}

Multiply \eqref{zcda3} by the integrating factor $e^{-\cl t}$ 
to obtain
\begin{equation}\label{cda51}
	\frac{d}{dt}\big(e^{-\cl t/}\|w\|^2\big)
	\le 2e^{-\cl t}\mu (w,AJ(U-v)).
\end{equation}
Now, integrate over $[t_n,t_{n+1})$ to obtain
\begin{align*}
	\|w&(t_{n+1}^-)\|^2 
	\le e^{\cl\delta}\|w_n\|^2
		+2\mu\int_{t_n}^{t_{n+1}} e^{\cl (t_{n+1}-t)}
			\big(w(t),AJ(U(t)-v(t))\big)dt\cr
	&\le
	e^{\cl\delta}\big\{(1+\delta\mu\cj)\|w(t_n^-)\|^2
        -2\delta\mu\big(
            w(t_{n}^-),
            AJ(U(t_{n})-v(t_n))
            \big)
        +\delta^2\mu^2\ck\big\}\cr
		&\qquad+2\mu\int_{t_n}^{t_{n+1}}
			e^{\cl(t_{n+1}-t)}(w(t),AJ(U(t)-v(t)))dt\cr
	&\le
	e^{\cl\delta}\big\{(1+\cj\delta\mu)\|w(t_n^-)\|^2
        +\delta^2\mu^2\ck\big\}\cr
		&\qquad+2\mu\int_{t_n}^{t_{n+1}}\Big\{
			e^{\cl(t_{n+1}-t)}\big(AJ(U(t)-v(t)),w(t)\big)\cr
			&\qquad\qquad\qquad\qquad -e^{\cl\delta}\big(
            AJ(U(t_{n})-v(t_{n})),
            w(t_{n}^-)
			\big)\Big\}
		dt\cr
	&\le
	e^{\cl\delta}\big\{(1+\cj\delta\mu)\|w(t_n^-)\|^2
        +\delta^2\mu^2\ck+2\mu(I_1+I_2+I_3+I_4)\big\}
\end{align*}
where 
$$	
	I_1=\int_{t_n}^{t_{n+1}} 
		\big(e^{\cl(t_n-t)}-1\big)
		\big(A^{1/2}J(U(t)-v(t)),A^{1/2} w(t)\big)dt,
$$
$$
	I_2=\int_{t_n}^{t_{n+1}}
		(AJ(U(t)-v(t)),w(t)-w(t_n^-))dt,
$$
$$
	I_3=\int_{t_n}^{t_{n+1}}
		\big(A^{1/2}J(U(t)-U(t_n)),
			A^{1/2}w(t_n^-)\big)dt
$$
and
$$
	I_4=-\int_{t_n}^{t_{n+1}}
		\big(A^{1/2}J(v(t)-v(t_n)),
			A^{1/2}w(t_n^-))dt.
$$
Let us estimate $I_1$ first.  Using \eqref{Jrho} and \eqref{JR}
followed by $(e^{\cl\delta}-1)\le \cl\delta e^{\cl\delta}$
yields
\begin{align*}
	I_1
	&\le\int_{t_n}^{t_{n+1}} 
		\big(1-e^{\cl(t_n-t)}\big)
		\big\|J(U(t)-v(t))\big\|\|w(t)\|dt\cr
	&\le (\rho_K+R_K)M_1 e^{-\cl\delta}
		\int_{t_n}^{t_{n+1}}
		\big(e^{\cl\delta}-1\big)dt
	\le \delta^2\cm
\end{align*}
where $\cm=(\rho_K+R_K)M_1 \cl$.

Second estimate $I_2$.  Using $J=P_\lambda J$ and
\eqref{psmooth} with $\alpha=2$ 
followed by Corollary~\ref{intdwbound} yields
\begin{align*}
	I_2
	&\le
		\int_{t_n}^{t_{n+1}}\big|AJ(U(t)-v(t))\big||w(t)-w(t_n^-)|dt\cr
	&\le
		\int_{t_n}^{t_{n+1}}\big|AP_\lambda J(U(t)-v(t))
			\big||w(t)-w(t_n^-)|dt\cr
	&\le
		\int_{t_n}^{t_{n+1}}\lambda^{1/2}
		\big\|J(U(t)-v(t))\big\||w(t)-w(t_n^-)|dt\cr
	&\le \lambda^{1/2}(\rho_K+R_K) \int_{t_n}^{t_{n+1}}|w(t)-w(t_n^-)|dt
	\le \delta^{3/2}\cn
\end{align*}
where $\cn=\lambda^{1/2}(\rho_K+R_K)M_E$.

Finally estimate $I_3$ and $I_4$
using Lemma~\ref{intdvbound} as
\begin{align*}
	I_3&\le\int_{t_n}^{t_{n+1}}
		\big\|J(U(t)-U(t_n))\big\|\|w(t_n^-)\|dt\cr
	&\le M_1\int_{t_n}^{t_{n+1}}
		\big\|J(U(t)-U(t_n))\big\|dt\le \delta^{3/2}\co
\end{align*}
where $\co=M_1 \rho_D$ and similarly
\begin{align*}
	I_4&\le\int_{t_n}^{t_{n+1}}
		\big\|J(v(t)-v(t_n))\big\|\|w(t_n^-)\|dt\cr
	&\le M_1\int_{t_n}^{t_{n+1}}
		\big\|J(v(t)-v(t_n))\big\|dt\le \delta^{3/2}\cp
\end{align*}
where $\cp=M_1R_D$.

It follows that
$$
	\|w(t_{n+1}^-)\|^2
	\le e^{\delta\cl}\big\{(1+\delta\mu\cj)\|w(t_n^-)\|^2
		+\delta^{3/2} \cq \big\}
$$
where $\cq=\dmax^{1/2}(\mu^2\ck+2\mu\cm)+2\mu(\cn+\co+\cp)$.

Setting $y_n=\|w(t_n^-)\|^2$ gives the inequality
$y_{n+1}\le \alpha y_n+\beta$
where $\alpha=e^{\delta\cl}(1+\delta\mu\cj)$ and
$\beta=\delta^{3/2} e^{\delta\cl}\cq $.  
Note that 
$$
	\alpha-1\ge\delta(\cl+\mu e^{\delta\cl}\cj)\wwords{and}
	\alpha^n\le e^{(t_n-t_0)(\cl+\mu\cj)}.
$$
Therefore, by induction
\begin{align*}
	\|w(t_n^-)\|
	&\le e^{(t_n-t_0)(\cl+\mu\cj)} \|w(t_0^-)\|
		+{e^{(t_n-t_0)(\cl+\mu\cj)}-1\over
		\cl+\mu e^{\delta\cl}\cj}\delta^{1/2}e^{\delta\cl}\cq\cr
	&\le e^{2\dmax(\cl+\mu\cj)} \|w(t_0^-)\|
		+{e^{2\dmax(\cl+\mu\cj)}-1\over
		\cl+\mu \cj}\delta^{1/2}e^{\dmax\cl}\cq.
\end{align*}
Since $\|w(t_0^-)\|=0$ it follows there is $\cs$ independent of
$\delta$ such that
$$
	\|w(t_n^-)\|\le \delta^{1/2}\cs
		\wwords{for} n= 1,\ldots,N.
$$

Again from \eqref{wnest} we have
$$
	\|w_n\|\le 
    (1+\delta\mu\cd)\delta^{1/2}\cs+\delta\mu\ce
	\le \delta^{1/2} C_{19}
$$
where $C_{19}=(1+\dmax\mu\cd)\cs+\dmax^{1/2}\mu\ce$.
Finally, \eqref{integrated} implies
$$
	\|w(t)\|^2
	\le e^{\cb \delta}\delta C_{19}^2+{\cc\over \cb}\big(e^{\cb\delta}-1\big)
	\le e^{\cb \delta}\delta C_{19}^2+
		\cc\delta e^{\cb\delta}\le \delta C_{20}
$$
where
$C_{20}=e^{\cb \dmax}(C_{19}^2+\cc)$.  Therefore
$\|w(t)\|\to 0$ as $\delta\to 0$ uniformly on $[t_0,T]$.
\end{proof}

\section{Concluding Remarks}\label{conclusion}

In the context of the incompressible two-dimensional 
Navier--Stokes equations
we have introduced a relaxation parameter $\kappa$ 
depending on the time 
$\delta$ between successive observations into the update step
of the spectrally-filtered direct-insertion method.
This
overcomes the lack of
synchronization when the observations are inserted too frequently.
The relationship between $\kappa$ and $\delta$ may be obtained numerically 
by varying these parameters 
independently and then minimizing the average error about an
ensemble of reference solutions.
It turns out that the relation $\kappa=\mu\delta$ for a particular
fixed $\mu$ works well over a wide range.

Analytically we have shown as $\delta\to0$
that the approximation corresponding
to the $\kappa$-relaxed discrete-in-time method converges to
the solution obtained by continuous nudging
on any finite time interval.
Although this does not 
imply $\|u(t)-U(t)\|\to 0$
as $t\to\infty$ for all sufficiently small $\delta$, the numerics 
suggest for suitable $h$, $\lambda$ and $\mu$ that $u$
does indeed synchronize with $U$.
A complete analysis proving such a result would be of great 
interest.

Data assimilation for the two-dimensional 
magnetohydrodynamic equations 
studied by Biswas, Hudson, Larios and Pie \cite{BJAY2017} and
Hudson and Jolly \cite{HJ2019}
provide an opportunity to further test the scaling 
$\kappa=\delta\mu$.
Let ${\bf S}$ be the semigroup 
corresponding to
the evolution of the two-dimensional magnetohydrodynamic equations
and ${\bf X}=(U_1,U_2,B_1,B_2)$ be the state of the velocity 
and magnetic fields.
Given the free-running solution ${\bf X}(t)={\bf S}(t){\bf X}_0$
consider interpolated observations ${\bf J}{\bf X}(t_n)$ 
at times $t_n$ 
taken 
on one component of each field.  Thus,
$$	
	{\bf J}(u_1,u_2,b_1,b_2)=(Ju_1,0,Jb_1,0)
$$
where $Ju_1$ and $Jb_1$ could be projections onto the $N$ lowest Fourier
modes or interpolants of the form $P_\lambda P_H I_h$ as in \eqref{interp}.
The analysis step corresponding to \eqref{relaxed} then becomes
$$
	{\bf x}_{n+1}
	=(I-\delta\mu {\bf J}) {\bf S}(t_{n+1}-t_n)
	{\bf x}_n
	+\delta\mu {\bf J} {\bf X}(t_n).
$$
Our conjecture is that ${\bf x}_n$ will as $\delta\to 0$ 
converge to the approximating solutions obtained by the 
continuous-in-time nudging method
given as Algorithm 2.3 in \cite{HJ2019}.

In the more general case where
$t_{n+1}-t_n=\delta_n$ and $\max_n \delta_n$ is small
\eqref{heursitica} reads as
$$
	\int_{t_n}^{t_{n+1}} \mu J(U-v) dt
		\approx \delta_n \mu J\big(U(t_n)-v(t_n)\big)
$$
which immediately suggests the time-dependent relaxation parameter 
$\kappa_n=\delta_n\mu$.
We expect similar results to those obtained here
also hold when \eqref{relaxed} 
is replaced by
$$
	u_{n+1}=(I-\delta_n\mu J)S(t_{n+1}-t_n)u_n +\delta_n\mu J U(t_{n+1}).
$$

We remark that the time between successive observations may not be 
uniform when there are missing observations.  Thus, it may 
happen that $\max_n \delta_n$ is not small but for most of the time
$\delta_n$ is small.
In this case taking
$\kappa_n=\min(1,\delta_n\mu)$ would avoid overrelaxation 
where $\kappa_n>1$.
This suggests another direction for future work:  Find
conditions on sequences of $\delta_n$ which for some values of 
$n$ are large and for which the approximation
$u$ numerically synchronizes with the reference solution $U$.

\appendix
\section{Appendix}\label{Tapriori}

In this appendix we follow the methods used to prove Theorem 5 in
\cite{AOT2014} to obtain explicit bounds on $v$ given by equation
\eqref{cda} that avoid imposing additional conditions on $h$ and
$\mu$ but hold only for a finite time interval $[t_0,T]$.
These bounds depend on $T$ and are valid 
for the case when $v$ does not synchronize with $U$.

We begin with the following existence theorem:

\begin{theorem}\label{vexist}
Suppose $J=P_\lambda P_H I_h$ where $I_h$ satisfies \eqref{interp}.
Then for any $h>0$ and $\mu>0$ the equations~\eqref{cda} with $f\in V$
and $v(t_0)\in V$ have unique 
strong solutions that satisfy
\beq\label{A11}
v\in C\big([t_0,T];V\big)\cap L^2\big((t_0,T);{\cal D}(A)\big)
\eeq
and 
\beq\label{A12}
\frac{dv}{dt} \in L^2\big((t_0,T);H\big)
\eeq
for any $T>t_0$.
\end{theorem}

The following analysis makes use of the Aubin Compactness Theorem
which we state here for reference as
\begin{theorem}\label{aubin}
Consider three separable reflexive Banach spaces 
$X_1\subset X_0\subset X_{-1}$ in which the inclusion 
$X_1\subset X_0$ is compact and the inclusion
$X_0\subset X_{-1}$ is continuous.  If $v_m$ is a sequence
such that
$$
    \big\|v_m\big\|_{L^2((t_0,T),X_1)}
\wwords{and}
    \Big\|{dv_m\over dt}\Big\|_{L^2((t_0,T),X_{-1})}
$$
are uniformly bounded in $m$.
There there is a subsequence $m_j$ and 
$v\in L^2\big((t_0,T),X_0\big)$ such that
$$
    \int_{t_0}^T \big\|v(\tau)-v_{m_j}(\tau)\big\|_{X_0}^2 d\tau\to 0
\wwords{as} j\to\infty.
$$
\end{theorem}
\noindent
Details are in Aubin \cite{Aubin1963}, the proof
of Lemma~8.2 in \cite{Constantin1988} or page 224
of \cite{Foias2001}.

Our analysis also makes use of a particular 
case of the Lions--Magenes 
theorem stated as Theorem 3 in Section 5.9.2 of  \cite{EvansPDE} or
Lemma 1.2 in Chapter III of \cite{Temam1977}. We state this result here
for reference as
\begin{theorem}\label{T11}
Suppose $X_1\subset X_0\subset X_{-1}$ is a Gelfand triple.
Thus, $X_1$ is densely and continuously embedded in $X_0$ and
$X_{-1}$ is the dual of $X_1$ with a pairing that extends 
the inner product on $X_0$.
Suppose $u\in L^2\big((t_0,T);X_1\big)$ 
with $du/dt\in L^2\big((t_0,T);X_{-1}\big)$. 
Then, after redefinition on a set of measure zero, 
the mapping $$t\to \|u(t)\|_{X_0}^2$$ is absolutely continuous with 
\beq
\frac{d}{dt}\|u(t)\|_{X_0}^2
	=2\Big\langle\frac{du(t)}{dt},u(t)\Big\rangle
\eeq
for almost every $t_0\le t\le T$. 
\end{theorem}

We turn now to the

\begin{proof}[Proof of Theorem \ref{vexist}]
Let $g=f+\mu JU$.  By Proposition~\ref{JUbound} 
\begin{equation}\label{gbound}
|g|\le |f|+\mu |JU| \le G_0\wwords{where} 
G_0=|f|+\mu\rho_J.
\end{equation}
Now proceed by the Galerkin method.

Let $P_m$ be the projection onto the $m$ lowest Fourier modes.
Note the projection $P_m$ employed in this proof is related 
to the spectral filter $P_\lambda$
introduced in \eqref{filter} such that $P_m=P_\lambda$ for
$m={\rm card}\{\, k\in\J : 0<|k|^2\le\lambda\,\}$.
Although the Galerkin method could be carried out
using $P_\lambda$,
for consistency with \cite{Temam1977}, \cite{Constantin1988},
\cite{Foias2001} and \cite{Robinson2001} we shall employ $P_m$.

Choose $m$ large enough that $P_mP_\lambda=P_\lambda$.
Thus, $P_m J=J$.
Now, consider the truncation of \eqref{cda} given by
\begin{equation}\label{galv}
\frac{dv_m}{dt}+\nu Av_m+P_m B(v_m,v_m)
	=P_m g-\mu Jv_m\\
\end{equation}
where $v_m(t_0)=P_m v_0$ and $t>t_0$.
Note $v_m\in P_m H$ is finite dimensional with dimension $m$.

Our goal is to find bounds on $v_m$  which 
are uniform in $m$ and hence to let $m\to\infty$ to obtain a 
solution for \eqref{cda}.
Taking inner product of \eqref{galv} with $v_m$ we have
\beq\label{cda2}
	\frac{1}{2}\frac{d}{dt}|v_m|^2+\nu \|v_m\|^2
	=(g,v)-\mu (Jv_m,v_m)\le \frac{1}{2\mu}|g|^2
		+\frac{\mu}{2}|v_m|^2-\mu (Jv_m,v_m).
\eeq
Estimate the last term of \eqref{cda2} using \eqref{Jinterp} as
\beq\label{cda2a}
\begin{aligned}
-\mu (Jv_m,v_m)&
	=\mu (v_m-Jv_m,v_m)-\mu|v_m|^2\le\mu|v_m-Jv_m||v_m|-\mu|v_m|^2\\
	&\le\frac{\mu}{2}\frac{\nu}{(\lambda^{-1}+c_1h^2)\mu}|v_m-Jv_m|^2
		+\frac{\mu}{2}\frac{(\lambda^{-1}+c_1h^2)\mu}{\nu}|v_m|^2-\mu|v_m|^2\\
	&=\frac{\nu}{2}\|v_m\|^2
		+\left(\frac{(\lambda^{-1}+c_1h^2)\mu^2}{2\nu}-\mu\right)|v_m|^2.
\end{aligned}
\eeq
Substituting \eqref{cda2a} back into \eqref{cda2} and multiplying by two 
then yields
\beq\label{cda3}
	\frac{d}{dt}|v_m|^2+\nu \|v_m\|^2	
		\le \frac{1}{\mu}|g|^2
			+\gamma |v_m|^2
\eeq
where $\gamma\ge(\lambda^{-1}+c_1h^2)\mu^2/\nu-\mu$.
We further assume $\gamma>0$.  
Note the case when $\gamma$ may be chosen negative 
was treated in \cite{AOT2014} and
leads to the synchronization of $v$ with $U$.

Now, drop the viscosity term $\nu \|v_m\|^2$ 
from \eqref{cda3} and apply \eqref{gbound}
to obtain
$$
\frac{d}{dt}|v_m|^2-\gamma|v_m|^2\le \frac{G_0^2}{\mu}
\wwords{or equivalently}
{d\over dt}\big(|v_m|^2e^{-\gamma t}\big)\le{G_0^2\over\mu}e^{-\gamma t}.
$$
Since $|v_m(t_0)|\le |v(t_0)|$, integrating from
$t_0$ to $t$ yields
\beq\label{cda8}
	|v_m(t)|^2\le |v(t_0)|^2 e^{\gamma (t-t_0)}+ \frac{G_0^2}{\gamma\mu}
		(e^{\gamma(t-t_0)}-1).
\eeq
Moreover, $\gamma>0$ implies $e^{\gamma t}\le e^{\gamma T}$ 
for $t\in[t_0,T]$.  It follows that
\beq\label{cda9}
\begin{aligned}
	|v_m(t)|
	&\le R_0 \wwords{where}
	R_0^2=|v(t_0)|^2e^{\gamma (T-t_0)}
		+\frac{G_0^2}{\gamma\mu}(e^{\gamma(T-t_0)}-1).
\end{aligned}
\eeq
We emphasize that $R_0$ depends on $T$ but 
holds uniformly in $m$.

Next, integrate \eqref{cda3} directly as
$$
	|v_m(T)|^2-|v(t_0)|^2+\nu \int_{t_0}^T\|v_m(\tau)\|^2d\tau
	\le \frac{1}{\mu}\int_{t_0}^T|g(\tau)|^2d\tau
		+\gamma\int_{t_0}^T|v_m|^2d\tau.
$$
From this \eqref{gbound} and \eqref{cda9} imply
\beq\label{cda11}
\begin{aligned}
\nu \int_{t_0}^T\|v_m(\tau)\|^2d\tau
	&\le \left(\frac{G_0^2}{\mu}+\gamma R_0^2\right)(T-t_0)+|v(t_0)|^2.
\end{aligned}
\eeq
It follows that
\beq\label{cda122a}
	\Big(\int_{t_0}^T\|v_m(\tau)\|^2d\tau\Big)^{1/2}\le \widetilde R_1
\words{where}
	\widetilde R_1^2=
		{1\over \nu}\left(\frac{G_0^2}{\mu}
		+\gamma R_0^2\right)(T-t_0)+\frac{1}{\nu}R_0^2.
\eeq

To estimate $\|v_m(t)\|^2$ and 
$\int_{t_0}^T|Av_m(\tau)|^2d\tau$ take inner product of 
\eqref{galv} with $Av_m$.  Thus,
\beq\label{cda13}
	\frac{1}{2}\frac{d}{dt}\|v_m\|^2
	+\nu |Av_m|^2+\big(B(v_m,v_m),Av_m\big)=(g,Av_m)-\mu (Jv_m,Av_m).
\eeq
From \eqref{bbthree} followed by Young's inequality with $p=4$ and 
$q=4/3$ we have
\beq\label{cda14}
\begin{aligned}
|(B&(v_m,v_m),Av_m)|
	\le c|v_m|^{1/2}\|v_m\|^{1/2}\|v_m\|^{1/2}|Av_m|^{1/2}|Av_m|\\
	&=\Big({6^{3/4}\over \nu^{3/4}}c|v_m|^{1/2}\|v_m\|\Big)
		\Big({\nu^{3/4}\over 6^{3/4}}|Av_m|^{3/2}\Big)
	\le \frac{54}{\nu^3}c^4|v_m|^{2}\|v_m\|^4+\frac{\nu}{8}|Av_m|^2.
\end{aligned}
\eeq
We also have
\beq\label{cda15}
	|(g,Av_m)|\le|g||Av_m|\le\frac{2}{\nu}G_0^2+\frac{\nu}{8}|Av_m|^2
\eeq
and recalling the definition of $\gamma$ that
\beq\label{cda16}
\begin{aligned}
	-\mu (Jv_m&,Av_m)=\mu (v_m-Jv_m,Av_m)-\mu\|v_m\|^2\\
	&\le\frac{\mu^2}{\nu}|v_m-Jv_m|^2+\frac{\nu}{4}|Av_m|^2-\mu\|v_m\|^2
	\le \gamma\|v_m\|^2 +\frac{\nu}{4}|Av_m|^2.
\end{aligned}
\eeq

Substituting \eqref{cda14}, \eqref{cda15} and \eqref{cda16}
into \eqref{cda13} yields
\begin{equation}\label{Aineq}
\frac{d}{dt}\|v_m\|^2+\nu |Av_m|^2
	\le\frac{108}{\nu^3}c^4|v_m|^{2}\|v_m\|^4
	+\frac{4G_0^2}{\nu}+2\gamma\|v_m\|^2
\end{equation}
and consequently
\beq\label{cda18}
	\frac{d}{dt}\|v_m\|^2
		-\Big(\frac{108}{\nu^3}c^4|v_m|^{2}\|v_m\|^2
		+2\gamma\Big)\|v_m\|^2\le\frac{4G_0^2}{\nu}
\eeq
for all $t\in[t_0,T]$. Define 
$$
	\psi_m(t)=\exp\left\{-\int_{t_0}^t 
		\Big(\frac{108}{\nu^3}c^4|v_m(\tau)|^{2}\|v(\tau)\|^2+2\gamma\Big)
		d\tau \right\}.
$$
Note that
\begin{align*}
	\int_{t_0}^t&
        \Big(\frac{108}{\nu^3}c^4|v_m(\tau)|^{2}\|v_m(\tau)\|^2+2\gamma\Big)
        d\tau \\
	&\le \frac{108}{\nu^3}c^4R_0^2\int_{t_0}^T 	
		\|v_m(\tau)\|^2d\tau+2(T-t_0)\gamma
	\le \frac{108}{\nu^3}c^4R_0^2\widetilde R_1^2+2\gamma(T-t_0)<\infty
\end{align*}
and therefore
$$
	\psi_m(t)\ge 
		\exp\left\{-\Big(\frac{108}{\nu^3}c^4R_0^2\widetilde R_1^2
			+2\gamma(T-t_0)\Big)\right\}
	\wwords{for all} t\in[0,T].
$$
Multiplying \eqref{cda18} by $\psi_m(t)$, integrating and
using that $\psi_m(t_0)=1$ and $\psi_m(t_0)\le 1$ yields
\beq\label{cda20}
\begin{aligned}
	\psi_m(t)\|v_m(t)\|^2-\|v_m(t_0)\|^2
		&\le\frac{4G_0^2}{\nu}\int_0^t\psi_m(\tau)d\tau
		\le\frac{4G_0^2}{\nu}(T-t_0).
\end{aligned}
\eeq
Since $\psi_m$ is decreasing on $[t_0,T]$,
it follows that
$$
	\|v_m(t)\|\le R_1
\wwords{where} R_1^2=\frac{1}{\psi_m(T)}\left\{\|v(t_0)\|^2
		+\frac{4G_0^2}{\nu}(T-t_0)\right\}.
$$

On the other hand, directly integrating inequality \eqref{Aineq} gives us
\beq\label{cda21}
\begin{aligned}
	\|v_m(T)&\|^2-\|v_m(t_0)\|^2+\nu \int_{t_0}^T|Av_m(\tau)|^2d\tau\\
		&\le\frac{108}{\nu^3}c^4\int_0^T
		\Big(|v_m(\tau)|^{2}\|v_m(\tau)\|^2+2\gamma\Big)\|v_m(\tau)\|^2 d\tau
			+ {4G_0^2\over\nu}(T-t_0).
\end{aligned}
\eeq
Estimate
\begin{align*}
	\int_{t_0}^T \Big(|v_m(\tau)|^{2}&\|v_m(\tau)\|^2
		+2\gamma\Big)\|v_m(\tau)\|^2 d\tau\cr
&\le
	(R_0^2R_1^2+2\gamma) \int_{t_0}^T \|v_m(\tau)\|^2 d\tau
\le
	(R_0^2R_1^2+2\gamma) \widetilde R_1^2.
\end{align*}
Substituting this into \eqref{cda21} yields
$$
	\Big(\int_{t_0}^T|Av_m(\tau)|^2d\tau\Big)^{1/2}\le \widetilde R_2
$$
where
$$
	\widetilde R_2^2=
		\frac{108}{\nu^4}c^4 (R_0^2R_1^2+2\gamma) \widetilde R_1^2
		+{4G_0^2\over\nu^2}(T-t_0)+{1\over\nu}R_1^2.
$$

Next we obtain uniform bounds on
$dv_m/dt$.
Setting $R_J=(\lambda^{-1}+c_1h^2)^{1/2}R_1+R_0$ 
as in the proof of Proposition \ref{JUbound} 
we obtain
that $|Jv_m|\le R_J$.
From \eqref{galv} and \eqref{bbone} holds
\begin{align*}
\Big|\frac{dv_m}{dt}\Big|
	&\le\nu|Av_m|+c|v_m|^{1/2}\|v_m\||Av_m|^{1/2}+G_0+\mu R_J\\
	&\le 2\nu|Av_m|+\frac{c^2}{4\nu}R_0 R_1^2+G_0+\mu R_J.
\end{align*}
Therefore,
$$
\Big|\frac{dv_m}{dt}\Big|^2
	\le 8\nu^2|Av_m|^2+2\Big(\frac{c^2}{4\nu}R_0 R_1^2+G_0+\mu R_J\Big)^2
$$
and we obtain
$$
	\Big(\int_{t_0}^T \Big|\frac{dv_m}{dt}\Big|^2\Big)^{1/2}
		\le R_T 
\words{where}
		R_T^2=8\nu^2 \widetilde R_2+
		2\Big(\frac{c^2}{4\nu}R_0 R_1^2+G_0+\mu R_J\Big)^2(T-t_0).
$$
Hence $dv_m/dt$ is bounded in $L^2((t_0,T),H)$ uniformly in $m$.

At this point we have shown that the solutions $u_m$ 
satisfy the bounds
$$
	\big\|v_m\big\|_{L^2((t_0,T),{\cal D}(A))} \le \widetilde R_2
\wwords{and}
	\Big\|{dv_m\over dt}\Big\|_{L^2((t_0,T),H)} \le R_T
$$
uniformly in $m$.
Since the inclusion ${\cal D}(A)\subseteq V$ is compact and
$V\subseteq H$ is continuous, then
by Theorem \ref{aubin}
there is a subsequence $m_j$ and $v\in L^2((t_0,T),V)$
such that 
$$
	\int_{t_0}^T \big\|v(\tau)-v_{m_j}(\tau)\big\|^2 d\tau\to 0
\wwords{as} j\to\infty.
$$
Taking additional subsequences if
necessary we may further suppose $v\in L^\infty((t_0,T);V)$
and that the weak time derivative $dv/dt\in L^2((t_0,T);H)$.

It follows that 
\begin{equation}\label{abnd1}
	v\in L^\infty((t_0,T);V)\cap L^2((t_0,T);{\cal D}(A))
	\wwords{and}
	dv/dt\in L^2((t_0,T);H).
\end{equation}
Moreover, since $J\colon V\to V$ is continuous then
for every $\phi\in H$ holds
\begin{equation}\label{cdastrong}
	\Big({dv\over dt},\phi\Big)+\nu(Av,\phi)+(B(v,v),\phi)
		=(g,\phi)-\mu (Jv,\phi)
\end{equation}
for almost every $t\in[t_0,T]$.

It remains to show such solutions are unique, contained
in $C([t_0,T];V)$ and depend continuously
on the initial data.
Let $v_1$ and $v_2$ be solutions satisfying \eqref{abnd1}
and \eqref{cdastrong}. Choose $K$ large enough such
that $\|v_1\|\le K$ and $\|v_2\|\le K$ for almost every $t\in[t_0,T]$.
Let $\tilde v=v_1-v_2$. Then $\tilde v$ satisfies

\beq\label{JV}
\frac{d\tilde v}{dt}+\nu A\tilde v+B(v_1,\tilde v)
	+B(\tilde v, v_2)=-\mu J\tilde v.
\eeq

By the Riesz representation theorem every element $\varphi^*\in V_2^*$
can be represented by $\varphi\in V_2$ 
such that 
$\langle\varphi^*,w\rangle=(\!(w,\varphi)\!)_2$ for all $w\in V_2$.
Setting $\psi_k=|k|^2\varphi_k$ we obtain $\psi\in H$ where
$$
	(\!(w,\psi)\!)
	=L^2\sum_{k\in\J} |k|^2 w_k\overline{\psi_k}
	=L^2\sum_{k\in\J} |k|^4 w_k\overline{\varphi_k}
	=(\!(w,\varphi)\!)_2
	=\langle\varphi^*,w\rangle.
$$
This identifies the dual of $V_2$ with $H$ in
such a way that extends the inner product of $V$.
Thus, ${\cal D}(A)\subset V\subset H$ is a Gelfand triple.

Now take inner product of \eqref{JV} with $A\tilde v$
and apply Theorem \ref{T11} to obtain
\beq\label{JV1}
	\frac{1}{2}\frac{d}{dt}\|\tilde v\|^2+\nu |A\tilde v|^2
		+(B(v_1,\tilde v),A\tilde v)+(B(\tilde v, v_2),A\tilde v)
		=-\mu (J\tilde v,A\tilde v).
\eeq
Estimate
$$
	|(B(v_1,\tilde v),A\tilde v)|
	\le c|v_1|^{1/2}\|v_1\|^{1/2}\|\tilde v\|^{1/2}|A\tilde v|^{3/2}
	\le \frac{27c^4}{4\nu^3\lambda_1}\|v_1\|^4\|\tilde v\|^2
		+\frac{\nu}{4}|A\tilde v|^2,
$$
$$
	|(B(\tilde v, v_2),A\tilde v)|
	\le c |\tilde v|^{1/2}\|v_2\||A\tilde v|^{3/2}
	\le \frac{27c^4}{4\nu^3\lambda_1} \|v_2\|^4  \|\tilde v\|^2
		+\frac{\nu}{4}|A\tilde v|^2
$$
and
$$
	-\mu (J\tilde v,A\tilde v)\le \mu |J\tilde v||A\tilde v|
		\le {\mu^2\over 2\nu}|J\tilde v|^2+{\nu\over 2}|A\tilde v|^2
		\le {\mu^2c_2^2\over 2\nu} \|\tilde v\|^2+{\nu\over 2}|A\tilde v|^2.
$$

Hence,
\begin{equation}\label{JV3}
	\frac{d}{dt}\|\tilde v\|^2
	\le \Big({27 c^4\over \nu^3\lambda_1}K^4+{\mu^2c_2^2\over\nu}\Big)
		\|\tilde v\|^2
\wwords{for all} t\in[t_0,T].
\end{equation}
Integrating \eqref{JV3} yields
$$
	\|\tilde v(t)\|^2
		\le \|\tilde v(t_0)\|^2 
		\exp\Big\{
			\Big({27 c^4\over \nu^3\lambda_1}K^4+{\mu^2c_2^2\over\nu}\Big)
		(t-t_0)\Big\}.
$$
Since $\tilde v(t_0)=v_1(t_0)-v_2(t_0)$
this shows strong solutions are unique and depend continuously 
on the initial data. 
\end{proof}

Note that the proof of
Lemma \ref{vexist} has already provided the bounds
$$
	\|v\|_\alpha\le R_\alpha \words{for} \alpha=0,1 \wwords{and}
\Big(\int_{t_0}^T \|v\|_\alpha^2\Big)^{1/2}\le \widetilde R_\alpha
	\words{for} \alpha=1,2.$$
To finish the proof of Theorem \ref{thmvbound} we need explicit
formulas for $R_2$ and $\widetilde R_3$.
To this end we further assume $v_0\in{\cal D}(A)$ and present the 
following formal estimates that 
if desired may be made 
rigorous using the same techniques as above.

Provided $f\in V$
and $U_0\in V$, it follows that $g=f+\mu J U\in V$ for all $t\ge 0$.
In particular,
$$
	\|g\|\le \|f\|+\mu \|JU\|\le G_1
\wwords{where}
	G_1=\|f\|+\mu \rho_K.
$$
Now take the inner product of equation \eqref{galv} with 
$A^2 v$ to obtain
$$
	{1\over 2}{d\over dt} |Av|^2+\nu |A^{3/2}v|^2+
		\big(A^{1/2}B(v,v),A^{3/2} v\big)
		=\big(A^{1/2}g,A^{3/2} v\big)
		-\mu \big(A^{1/2}J v,A^{3/2} v\big).
$$
Estimate using \eqref{bvnorm} as
\begin{align*}
	\big|\big(A^{1/2}&B(v,v),A^{3/2} v\big)\big|
	\le \|B(v,v)\| |A^{3/2}v|
	\le c |v|^{1/2}\|v\|^{1/2}|Av|^{1/2}|A^{3/2}v|^{3/2}\cr
	&\le {9^3\over 2^5}{c^4\over\nu^3}
		|v|^2\|v\|^2|Av|^2 +
		{\nu\over 6}|A^{3/2}v|^2
	\le {9^3\over 2^5}{c^4\over\nu^3}
		R_0^2 R_1^2 |Av|^2 +
		{\nu\over 6}|A^{3/2}v|^2.
\end{align*}
Also
$$
	\big|\big(A^{1/2}g,A^{3/2} v\big)\big|
		\le \|g\||A^{3/2}v|
		\le {3\over 2\nu} G_1^2+{\nu\over 6}|A^{3/2}v|^2
$$
and
$$
	\mu \big|\big(A^{1/2}J v,A^{3/2} v\big)\big|
		\le \mu\|Jv\||A^{3/2}v|
		\le {3\mu^2\over 2\nu}R_K^2 + {\nu\over 6}|A^{3/2}v|^2.
$$

It follows that
\begin{equation}\label{avdiffeq}
	{d\over dt}|Av|^2+\nu|A^{3/2}v|^2 \le c_4 |Av|^2+c_5
\end{equation}
where
$$
	c_4={9^3\over 2^4}{c^4\over\nu^3}
        R_0^2 R_1^2
\wwords{and}
	c_5={3\over \nu}\big(G_1^2+\mu^2R_K^2\big).
$$
Dropping the second term of the left, multiplying by $e^{-c_4t}$ 
and integrating over $[t_0,t)$ yields
$$
 	|Av|^2\le e^{c_4 (t-t_0)} |Av_0|^2
		+ {c_5\over c_4}\big(e^{c_4 (t-t_0)}-1\big).
$$
Therefore, for $t\in [t_0,T]$ holds
$$
	|Av|^2\le R_2
	\wwords{where}
	R_2=e^{c_4(T-t_0)}|Av_0|^2 + {c_5\over c_4}\big(e^{c_4(T-t_0)}-1\big).
$$
This is the $\alpha=2$ bound needed for \eqref{vbound1}.

Next, directly integrate \eqref{avdiffeq} over $[t_0,T]$ to obtain
$$
	|A v(T)|^2-|Av_0|^2+\nu \int_{t_0}^T |A^{3/2}v|^2\le
		c_4\int_{t_0}^{T} |Av|^2 + (T-t_0) c_5
$$
so that
$$
	\Big(\int_{t_0}^T
		|A^{3/2}v|^2\Big)^{1/2}\le \widetilde R_3
\wwords{where}
\widetilde R_3^2={1 \over \nu}\big(R_2^2 + c_4\widetilde R_2^2+(T-t_0)c_5\big).
$$
Finally, note that $[t_n,t_{n+1}]\subseteq [t_0,T]$ implies
$$
	\Big(\int_{t_n}^{t_{n+1}} \|v\|_\alpha^2\Big)^{1/2} 
		\le \Big(\int_{t_0}^T \|v\|_\alpha^2\Big)^{1/2}
	\le \widetilde R_\alpha.
$$
Therefore, the bounds $\widetilde R_\alpha$ for $\alpha=1,2,3$ obtained 
in this appendix satisfy \eqref{vbound2}.

\bibliographystyle{plain}

\end{document}